\documentclass{article}

\usepackage{lmodern,lscape,textcomp,url,placeins}
\usepackage[T1]{fontenc}

\usepackage{amsmath,amssymb,dsfont,units,booktabs}
\usepackage{amsthm}
\usepackage{graphicx, color}
\usepackage{color}
\usepackage[utf8]{inputenc} 
\usepackage{marginnote}
\usepackage{todonotes}
\usepackage{sidenotes}
\usepackage{appendix}

\newtheorem{theorem}{Theorem}

\newcommand{\vt}{\mathbf{v}}

\newcommand{\xt}{\mathbf{x}}

\newcommand{\nt}{\mathbf{n}}

\newcommand{\taut}{\pmb{\tau}}

\newcommand{\sigmat}{\pmb{\sigma}}
\newcommand{\phit}{\pmb{\phi}}
\newcommand{\epsilont}{\dot{\pmb{\epsilon}}}

\newcommand{\te}{\pmb{e_r}}
\newcommand{\R}{\color{black}}
\newcommand{\RR}{\color{black}}
\newcommand{\BB}{\color{black}}
\newcommand{\B}{\color{black}}

\newcommand{\G}{\color{black}}
\newcommand{\GG}{\color{black}}
\newcommand{\BR}{\color{black}}
\renewcommand{\div}{\operatorname{div}}
\newcommand{\tr}{\operatorname{tr}}

\newcommand{\rhs}{ \mathcal{R}}

\title{Sea-ice dynamics on triangular grids} 
\author{Carolin Mehlmann, Peter Korn}
\date{January 2021}
\begin{document}

\maketitle
\date
\begin{abstract}
 We present a stable discretization of sea-ice dynamics on triangular grids that can straightforwardly be coupled to an ocean model on a triangular grid with Arakawa C-type staggering. The approach is based on a nonconforming finite element framework, namely the Crouzeix-Raviart finite element. As the discretization of the viscous-plastic and elastic-viscous-plastic stress tensor with the Crouzeix-Raviart finite element {\RR produces oscillations in the velocity field, we introduce an edge-based stabilization. {\RR To show that the stabilized Crouzeix-Raviart} approximation is qualitative consistent with the solution of the continuous sea-ice equations, we derive} {\BR a $H^1$-estimate}.
 In a numerical analysis we show that  the stabilization is fundamental to achieve stable approximation of the sea-ice velocity field.  
\end{abstract}

\section{Introduction}
Sea-ice, located at high-latitudes and at the boundary between ocean and atmosphere, plays an important role in the climate system. 
Modelling the complex mechanical and thermodynamical behaviour of sea-ice at a broad range of spatio-temporal scales poses a manifold of challenges. 
Freezing sea water forms a composite of pure ice, liquid brine, air pockets and solid salt. The details of this formation depend on the laminar or turbulent environmental conditions. 
This composite responds differently to heating, pressure or mechanical forces than for example the (salt-free) glacial ice of the ice sheets. 
Climate models need to describe the dynamics of sea-ice on large scales and couple the large-scale sea-ice models to ocean general circulation models. This is also the perspective we pursue in this work.

This paper treats the problem of formulating the discrete sea-ice dynamics in a way such that the internal sea-ice dynamics are captured well while at the same time the external coupling to the ocean is accomplished in a natural way. The modelling problem we aim to solve consists in choosing approximation spaces that capture (compressible) sea-ice dynamics as well as the (incompressible) ocean dynamics and that allow a minimal-invasive coupling between the two models that avoids interpolations or projections. 
In global ocean modelling we observe a trend towards a Arakawa C-type staggering of variables, where scalar variables are described as piecewise constant functions that are located at the center of the grid cell and velocity fields are represented by normal components of the velocity vector along the cell boundary 
(see e.g. \cite{MOM6, KORN_OPE,NEMO, MPAS-O}). This trend reflects the increase in computational power towards high-resolution simulation where 
the efficient C-type staggering has advantageous discrete wave propagation properties once the Rossby radius is resolved over large part of the global domain.
For sea-ice dynamics the small stencil of C-type staggering has the advantage that transport in narrow straits, which are only one cell wide is possible {\R \cite{Boullion2009}}. The modelling challenge stems from the fact that sea-ice dynamics  
requires the full strain rate tensor $\frac{1}{2}(\nabla \vt+ \nabla \vt^T)$ which is difficult to discretize if only partial information about the velocity field is available.

We consider this problem on triangular grids in the context of the ocean general circulation model ICON-O \cite{KORN_OPE}.  
ICON-O uses a refined triangular mesh of an icosahedron that is approximating the surface of the sphere as described in \cite{KORN_OPE}. 
On triangular grids the C-type staggering is equivalent to the lowest order  Raviart-Thomas finite element (RT-0).
The space of the Raviart-Thomas finite element is not rich enough to approximate the full strain rate tensor \cite{Acosta2011}. 

{\R We enlarge} the approximation space by including the tangential velocity at the midpoint of an edge. This variable arrangement allows the desired natural coupling to the 
underlying ocean variables on the same grid. This enrichment of the Raviart-Thomas element results in the specification of the complete velocity vector at edge midpoints and equals to the first order nonconforming Crouzeix-Raviart finite element. The  Crouzeix-Raviart element is a classical finite element that has been applied for approximations of the {\RR Poisson} problem and the Navier-Stokes equations. {Lietear et al.~\cite{Lietaer2008} used the Crouzeix-Raviart  element to discretize the Canadian Arctic Archipelago in an uncoupled sea-ice model.}

{\RR 
{\BB A direct application of the Crouzeix-Raviart element to the sea-ice equations leads to an unstable approximation of the sea ice velocity.}
The reason is an  instability that has its origin in the discretization of the symmetric strain rate tensor in the sea-ice rheology. More precisely, the  
discretization of the strain rate tensor with the Crouzeix-Raviart element 
causes a non-trivial null space, in addition to the trivial null space, spanned by solid-body rotations. Mathematically this problem can be characterized by the fact 
that the Crouzeix-Raviart element violates the Korn inequality. This is true even if the trivial null space is removed from the function space under consideration. This is further discussed in Section \ref{sec:disc}.  

The way out of this dilemma comes from the observation that the Crouzeix-Raviart element fulfills a {\it generalized version of Korn's inequality}, which is shown in \cite{Brenner2004}. Motivated by Korn's inequality in its generalized form we add a stabilizing term to the sea-ice momentum equation {\GG in spirit of an interior penalty method}. A similar stabilization was derived by Hansbo and Larson \cite{HansboLarson2003}. {\BR To shorten the notation we refer to the Crouzeix-Raviart discretization of this stabilized momentum equation, as the {\it stabilized CR element}.}

{\BB We demonstrate through numerical experiments that with this {\RR edge} based stabilization method the Crouzeix-Raviart element is suitable to discretize the viscous-plastic and elastic-viscous-plastic sea-ice model}. {\BR To show that the stabilized Crouzeix-Raviart approximation is qualitatively consistent with the solution of the continuous model equations, we derive a weighted $H^1$-estimate} for the viscous-plastic and elastic-viscous-plastic sea-ice model and evaluate it numerically. We show that only the  stabilized form of the Crouzeix-Raviart element fulfills the weighted $H^1$-estimate}. 

All numerical experiments are carried out in the framework of the ocean general circulation model ICON-O \cite{KORN_OPE} that {\RR operates} on a triangular C-grid. 


The paper is structured as follows. In Section 2 we introduce the sea-ice model in a strong and variational form and derive the weighted 
$H^1$-estimate for the viscous-plastic and elastic-viscous-plastic sea-ice momentum equations. In Section 3 we introduce the realization of  Crouzeix-Raviart element in ICON-O and describe the stabilization for the viscous-plastic and elastic-viscous-plastic model. In Section 4 we numerically analyze and validate the stabilized Crouzeix-Raviart element. {\BR In Section 5 the paper ends with a conclusion.}

\section{Model description}
 The motion of sea-ice is prescribed in a two-dimensional framework \cite{Lepp2010}.
The momentum of sea-ice is modelled as
\begin{align}\label{eq:mo}
\rho h \partial_t \vt= \div(\sigmat)+ F,
\end{align} 
where {\B $\vt$ is the horizontal sea ice velocity}, $\rho$ the sea-ice density, and $h$ the mean sea-ice thickness. All external forces are collected in $F$,
\begin{align*}
F=\tau-\rho h g \nabla H_d-\rho h f_c \te \times \vt,
\end{align*}
where $f_c$ is the Coriolis parameter, $g$ is the gravity, $\te$ the unit normal vector to the surface and $H_d$ is the ocean surface height. {\RR  We follow Coon \cite{Coon1980}} and use
\begin{align*}
 \rho h g \nabla H_d \approx - \rho h f\te \times \vt_{w},
\end{align*}
where $\vt_w$ is the ocean {velocity}.
The stresses due wind and ocean are {\BR given by}
\begin{align}
\tau&:=\rho_a C_a\|\vt_a\| \vt_a +\rho_w C_w \|\vt-\vt_w\| (\vt_w-\vt),\label{stress_def}
\end{align}
where $\rho_a, \rho_w$ are the air and water densities, $C_a,C_w$  are the wind and water drag coefficients and $\vt_a$ { is the atmospheric velocity}. 
To ensure that the formulation of the momentum equations is consistent with the formulation of free drift, $\tau$ should be multiplied with the sea-ice concentration $A$ \cite{Connolley2004}. We focus in this manuscript on the numerical development and do not take this scaling in our numerical experiments into account.

\paragraph{The viscous-plastic sea-ice rheology (VP)}
The internal stresses $\sigmat$ are related to the strain rate $\epsilont:= \frac{1}{2} (\nabla \vt+\nabla\vt^T)$ by the viscous-plastic rheology
\begin{align}\label{stress}
\sigmat = 2 \eta \epsilont +( \zeta- \eta) \tr(\epsilont) I -\frac{P}{2}I,
\end{align}
with the viscosities $\eta, \zeta$, and the ice strength $P$. 
Following Hibler the ice strength $P$ is modeled as
\begin{align*}
P=P^*h\exp(-C(1-A)),
\end{align*} 
with $P^*>0$ the ice strength parameter and $C>0$ the ice concentration parameter \cite{Hibler1979}.
In the viscous-plastic sea-ice model introduced by Hibler \cite{Hibler1979}, the viscosities are derived from a elliptic yield curve with eccentricity e=2 and a normal flow rule. They are modelled as 
\begin{align}\label{viscosities}
\zeta=\frac{P}{2\Delta}, \quad \eta=\zeta e^{-2}, \quad \Delta= \sqrt{\big(  2e^{-2} \epsilont':\epsilont'+\tr(\epsilont)^2\big)},
\end{align}
where we apply the decomposition of the strain rate tensor $\epsilont=\epsilont'+\frac{1}{2} \tr(\epsilont) I$ into the deviatoric part $\epsilont'$ and into its trace $\tr(\epsilont) I$.
As described by Hibler \cite{Hibler1979} the plastic viscosity {\RR is} limited by a viscous regime given as
{\GG
\begin{align*}
\zeta\leq \zeta_{max},\quad \zeta_{max}=\frac{P}{2\Delta_\text{min}}, \quad \Delta_\text{min}=2\cdot 10^{-9} \text{1/s}. 
\end{align*}
}
The limitation avoids that  $\zeta, \eta \to \infty$  for $\Delta \to 0$. To regularize the transition from the viscous to the plastic regime we follow Harder~\cite{Kreyscher2000} and use
\begin{align*}
\Delta=\sqrt{\Delta^2_\text{min}+ 2e^{-2} \epsilont':\epsilont'+\tr(\epsilont)^2.
}
\end{align*}

\paragraph{The elastic-viscous-plastic sea-ice rheology (EVP)}
The elastic-viscous-plastic model was introduced to regularize the VP rheology, 
such that the VP model results from the EVP model { for $\partial_t \sigmat=0$}. We reformulate the viscous-plastic model to \begin{align*}
\frac{1}{2\eta}\sigmat+\frac{\eta-\zeta}{4\eta\zeta}\tr(\sigmat)I+\frac{P}{4\zeta}I=\frac{\zeta}{T}\epsilont
\end{align*}
and add an artificial elastic strain  behaviour with a parameter E 
\begin{align*}
\frac{1}{E}\partial_t \sigmat +\frac{1}{2\eta}\sigmat +\frac{\eta-\zeta}{4\eta\zeta}\tr{\sigmat}I+\frac{P}{4\zeta}I=\epsilont.
\end{align*}
Hunke \cite{Hunke2001}  introduced $T_\text{EVP}=\frac{\zeta}{E}$ such that the elastic-viscous-plastic model {becomes }
\begin{align}\label{EVP}
\partial_t \sigmat +\frac{e^2}{2T_\text{EVP}}\sigmat +\frac{1-e^2}{4T_\text{EVP}}\tr{\sigmat}I+\frac{P}{4T_\text{EVP}}I=\frac{\zeta}{T_\text{EVP}} \epsilont.
\end{align}
The EVP model allows a fully explicit discretization in time with relatively large time steps such that it gains numerical efficiency compared to the VP model \cite{Hunke2001}.
In Section \ref{sec:disc} we introduce a modified EVP approach ({\bf mEVP}), a pseudo-time solver for the VP sea-ice model based on the EVP formulation.

All constants used in the momentum equation are {summarized} in Table \ref{Con}.
    \begin{table}[t]
  \begin{center}
    \begin{tabular}{|l|l|l|}
      \hline
      \text{Parameter} & \text{Definition} & \text{Value}\\
      \hline
      $\rho$ & sea-ice density &$\unit[900]{kg/m^{3}}$\\
      $ \rho_\text{a}$& air density & $\unit[1.3]{kg/m^{3}}$\\
      $\rho_\text{w}$& water density & $\unit[1026]{kg/m^{3}}$\\
      $C_\text{a}$&air drag coefficient &$\unit[1.2] \cdot {10^{-3}}$\\
      $C_\text{w}$&water drag coefficient& $\unit[5.5] \cdot {10^{-3}}$\\
      $f_c$ &Coriolis parameter&$\unit[1.46] \cdot \unit[10^{-4}]{s^{-1}}$ \\
      $P^{\star}$&ice strength parameter&$\unit[27.5]\cdot \unit[ 10^{3}]{N/m^2}$\\
      $C$&ice concentration parameter&$20$\\
      $e$&{\GG eccentricity of elliptic yield curve} &$2$\\
      \hline
    \end{tabular}
 \caption{Physical parameters of the momentum equation\label{Con}.}
  \end{center}
\end{table} 
\paragraph{Transport equations}
The mean sea-ice thickness and sea-ice concentration are advected in time by
{\B
\begin{align}\label{eq:transport}
\partial_t h+\div(\vt h)=Q_h, \quad \partial_t A+\div(\vt A)=Q_A,
\end{align}}
where $A$ is limited from above by $1.0$. Here we skip thermodynamic  source terms on the right hand side of the transport equations as our analysis focuses on the sea-ice dynamics and set $Q_h=0$ and $Q_A=0$. 

\paragraph{Weak formulation}
The weak form of the sea-ice equations is given as 
\begin{align}\label{eq:system}
(\rho h \partial_t \vt, \phit)-(F, \phit) +(\sigmat, \nabla \phit) &=0\\
(\partial_t h+\div(\vt h), \phi)&=0,\\
(\partial_t A+\div(\vt A), \phi)&=0,
\end{align}

where $(\cdot,\cdot)$ denotes the $L^2$- inner product on $\Omega$.

\subsection{{\BR Gradient} estimates for the viscous-plastic and the elastic-viscous-plastic sea-ice rheology}\label{sec:energy}
 In this section we {\BR derive an estimate for the weighted $H^1$-norm of velocity fields that satisfy the sea-ice momentum equation either with the viscous-plastic {or elastic-viscous-plastic rheology}.} Energy methods are an important mathematical tool for proving well-posedness and regularity results for nonlinear partial differential equations. For the equations of sea-ice dynamics questions such as existence and uniqueness of solutions is open. 
{\RR The estimates that we derive in this section are applied in our numerical experiments in Section \ref{sec:VPEVP} as a diagnostic to investigate {\BB if the CR approximation of  the stabilized momentum equation is qualitatively consistent with the solution of the  (elastic)-viscous-plastic sea-ice momentum equation.}

For the derivation of the weighted $H^1$-estimates we refer to Appendix \ref{app}.   
The assumptions we made in the derivation of the {\GG estimate} such as a fixed mean sea-ice thickness and sea-ice concentration, are incorporated in our experimental configuration in Section \ref{sec:VPEVP}.} 
{\B We assume a linearized water drag term (see \cite[Chapter 6.1.4]{Lepp2010}), {and use a
fixed scaling speed indicated with the subscript $0$}. Thus, the {\RR linear wind} drag is given by
\begin{align}\label{lin:drag}
C_{w1}=C_w| (\vt_w-\vt)|_0. 
\end{align}} {\RR {\BB The bound on the gradient of the sea ice velocity derived in Theorem 1 and Theorem 2 in the Appendix \ref{app} relies on the validity of Korn's inequality. The CR element violates the Korn inequality even if the trivial null space of the strain rate tensor, that is introduced by the solid body rotations, 
is removed. Therefore, we focus our numerical analysis on a weighted gradient of the velocity, defined by}
\begin{align}\label{eq:energy}
E(\vt):=\int_I \zeta_{\text{min}} \|\nabla \vt \|^2 dt, \quad {\GG \zeta_\text{min}:=\min_{(\xt,t) \in \Omega\times I} \zeta(\xt,t).}
\end{align}
Theorem 1 and Theorem 2 state that the weighted gradient $E(\vt)$ satisfies the following inequalities 
\begin{equation}\label{eq:est}
  \begin{aligned}
 \textbf{VP:} && \frac{c_k}{4} E(\vt) &\leq \int_I\frac{c_p}{c_k\zeta_\text{min}} 
  \|\rhs\|^2 dt 
+\|\rho h \vt(0) \|^2,\\
\textbf{EVP:} &&\frac{c_k}{2T_\text{EVP}}E(\vt) &\leq  \int_I  c \|  \partial_t \rhs+\rhs \|^2 dt\\
 &&  &+\frac{\zeta}{4}\|\epsilont(0)\|^2 + a \|\partial_t \vt(0)\|^2 +b \| \vt(0)\|^2,
\end{aligned}  
\end{equation}
}{\RR 
{where $\rhs:=\rho_a C_a\|\vt_a\| \vt_a+\rho_wC_{w1}\vt_w-\rho h f\te \times \vt_{w}$
and where $c_k$ and $c_p$ are positive constants that depend on the domain. In case of the EVP rheology the positive constant 
$a,b$ and $c$ are given as follows
\begin{align*}
    a=\frac{\rho h T_\text{EVP}}{4}+ \frac{1}{2\gamma}\rho h,\quad
    b= \rho_w C_{w1}  +\frac{1}{T_\text{EVP}}\rho h + \frac{\gamma}{2}\rho h \quad
    c=  \frac{c_p \zeta_\text{min}}{2 c_k T_\text{EVP}} + \frac{T_\text{EVP}}{4 \gamma_2}.
\end{align*}}
}
{\GG The estimates in (\ref{eq:est}) state that the weighted gradient of the solution of the continuous sea-ice equations is bounded by the data. Therefore, we expect the same for discretized solution and analyze it in Section \ref{sec:VPEVP}.}

\section{Discretization}\label{sec:disc}
In this section we describe the spatial and temporal discretization of the VP and EVP model. We introduce an edge-based stabilization for both models discretized with the Crouzeix-Raviart finite element. 
\paragraph{Time discretization}
To solve the coupled sea-ice system (\ref{eq:system}) it is standard to use a splitting approach in time. As described by Lemieux et al.~\cite{Lemieux2014}
we first  compute the solution of the sea-ice  momentum equation (\ref{eq:mo}), followed by the solution of the  transport equations (\ref{eq:transport}). Ip et al.~\cite{Hibler1991} pointed out that a fully explicit time
stepping scheme for the momentum equation with a VP rheology would
require a small time step of less than a second - even on a
grid resolution as coarse as  $\unit[100]{km}$.
Therefore the authors recommended an implicit treatment in time. An implicit discretization asks for implicit solution methods such as a Picard solver \cite{Hibler1997} or Newton like methods \cite{Lemieux2009, MehlmannRichter2016newton}. So far the applied solvers are difficult to parallelize as efficient linear solver are missing \cite{Losch2014}. To avoid an implicit discretization Hunke and Dukowicz \cite{Hunke2001, Hunke1997} introduced the EVP model, where they add an artificial elastic term to the VP rheology, to allow an explicit discretization of the momentum equation with relatively large time steps. However, the EVP model produces large differences compared to approximations of the VP model. Thus, Kimmritz et al.~\cite{Kimmritz2015} and Boullion et al.~\cite{Boullion2013} developed, based on the EVP model formulation, explicit pseudo-time stepping methods that converge against the solution of the VP model. We observe in numerical experiments that spatial discretization errors dominate the temporal discretization errors. Thus a first order time stepping methods are sufficient to discretize the sea-ice momentum equation \cite{Mehlmann2019, Mehlmann2020}.
A second-order time stepping scheme to solve the VP model is described in \cite{Lemieux2014}.

For the temporal discretization of the momentum equation (\ref{eq:mo}) we apply a first order semi-implicit time stepping scheme. We introduce the time partitioning $0=t_0< ...<t_N=T$ and the time step size $k:=t_{n}-t_{n-1}$. Let $ \vt^n:=\vt(t_n),\, h^n:=h(t_n),\, A^n:=A(t_n)$. 
\paragraph{VP model}
The discretized viscous-plastic momentum equation reads as
\begin{align}\label{eq:momVP}
\begin{aligned}
&\Big(\rho h^{n-1} \frac{\vt^n-\vt^{n-1}}{k}, \phit\Big)=
\Big(\rho h^{n-1} f_c e\times (\vt^{n-1}_w -\vt^{n-1}),\phit\Big)
- (\sigmat^{n},\nabla \phit)
+\tau^n_{VP}\\
&\text{with }\\
&\tau^n_{VP}:=\Big(\rho_a C_a\|\vt^{n-1}_\text{a} \| \vt^{n-1}_\text{a},\phit\Big )
+\Big (\rho_\text{w}C_w\| \vt^{n-1}_w-\vt^{n-1}\|( \vt^{n-1}_w-\vt^{n}), \phit\Big). \\
\end{aligned}
\end{align}

\paragraph{EVP model}
To solve the elastic-viscous-plastic model we sub-cycle the momentum equation. Let $t^{n-1} \leq t^{s-1}< t^s \leq t^n$ and $k_s:=t^s-t^{s-1}$. Then the subscycled momentum equation reads as
\begin{align}\label{eq:EVP}
\begin{aligned}
&\Big (\rho h^{n-1} \frac{\vt^s-\vt^{s-1}}{k_s},\phit \Big)=\Big( \rho h^{n-1} f e \times (\vt_w^{n-1}-\vt^{s-1}),\phit\Big)+\Big(\sigmat^s,\nabla \phit\Big)
+\tau^s_{EVP},\\
&\text{with }\\
&\tau^n_{EVP}:=\Big( \rho_a C_a\|\vt^{n-1}_\text{a} \| \vt^{n-1}_\text{a},\phit\Big)
+\rho_\text{w}C_w\| \vt^{n-1}_w-\vt^{s-1}\| (\vt^{n-1}_w-\vt^{s}),\phit \Big).
\end{aligned}
\end{align}
 The elastic-viscous-plastic stress is calculated via
\begin{align*}
\frac{\sigmat_{1}^s-\sigmat_{1}^{s-1}}{k_s}+\frac{\sigmat^s_{1}}{2T_\text{EVP}}&=\frac{\zeta^{s-1}(\epsilont^{s-1}_{11}+\epsilont^{s-1}_{22})}{T_\text{EVP}} -\frac{P^{n-1}}{2T_\text{EVP}},\\
\frac{\sigmat^s_{2}-\sigmat^{s-1}_{2}}{k_s}+\frac{4\sigmat^s_{2}}{2T_\text{EVP}}&=\frac{\zeta^{s-1}(\epsilont^{s-1}_{11}-\epsilont^{s-1}_{22})}{T_\text{EVP}},\\
\frac{\sigmat^s_{12}-\sigmat^{s-1}_{12}}{k_s}+\frac{4\sigmat^{s}_{12}}{2T_\text{EVP}}&=\frac{\zeta^{s-1}\epsilont^{s-1}_{12}}{T_\text{EVP}},
\end{align*}
with  $\sigmat_1=\sigmat_{11} +\sigmat_{22}$,  $\sigmat_{2}=\sigmat_{11}-\sigmat_{22}$, $\epsilont_1=\epsilont_{11}+\epsilont_{22}$, $\epsilont_{2}=\epsilont_{11}-\epsilont_{22}$ and $\zeta^{s-1}:= \zeta(\vt^{s-1}),\, \epsilont^{s-1}:=\epsilont(\vt^{s-1})$.
For all computations done in this paper we use the tuning parameter $T_\text{EVP}=100$ {\RR seconds}. The number of sub-cycles $s=1,...,N_\text{evp}$ is usually a large number around $100$ or more \cite{Hunke1997} and  $k_s=\frac{k}{N_\text{evp}}$, where $k$ is the larger time step of the advection.  

\paragraph{mEVP solver}
We approximate the viscous-plastic stress tensor with an elastic-viscous-plastic formulation. The elastic-viscous-plastic model is sub-cycled in time, such that the formulation converges in time against the viscous-plastic formulation. Each sub-iteration $s$ of the momentum equation reads as 
\begin{align}\label{eq:mEVP}
\begin{aligned}
&\Big(\beta\frac{\rho h^{n-1}}{k_s}(\vt^s-\vt^{s-1}),\phit\Big)
=\Big (\rho h^{n-1}\frac{-\vt^s+\vt^n}{k_s}, \phit \Big) 
- \Big (\sigmat^s,\nabla \phit \Big)\\
&\qquad\qquad\qquad\qquad\qquad - \Big(\rho h^{n-1} f_c \vec e\times (\vt^{s-1} -\vt^{n-1}_w),\phit \Big)+\tau^s_{mEVP}\\
&\text{with }\\
&\tau^s_{mEVP}:=  \Big( \rho_a C_a\|\vt^{n-1}_\text{a} \| \vt^{n-1}_\text{a}, \phit \Big)
+
\Big (\rho_\text{w}C_w\| \vt^{n-1}_w-\vt^{s-1}\|(\vt^{n-1}_w-\vt^{s}),\phit \Big),
\end{aligned}
\end{align}
where we time step the stress tensor $\sigmat^s$  as
\begin{align*}
\alpha (\sigmat_1^s-\sigmat_1^{s-1})&=\sigmat^{s-1}_1+2 \zeta^{s-1}(\epsilont^{s-1}_1-P^{n-1}),\\
\alpha(\sigmat_2^s-\sigmat_2^{s-1})&=\sigmat^{s-1} \frac{\zeta^{s-1}}{2}\epsilont^{s-1}_2,\\
\alpha(\sigmat_{12}^s-\sigmat_{12}^{s-1})&=\sigmat^{s-1} \frac{\zeta^{s-1}}{2}\epsilont^{s-1}_{12}.
\end{align*}
A more detailed description of the mEVP solver can be found in the work of Boullion et al.~\cite{Boullion2013}. Here $\alpha$ and $\beta$ are large constants. As analyzed by Boullion et al.~\cite{Boullion2013} and Kimmritz et al.~\cite{Kimmritz2015}, the product $\alpha\beta$ should be sufficient large to satisfy the CFL-criterion. If not further specified we follow Koldunov et al.~\cite{Koldunov} and use $\alpha=\beta=500$ in this paper.

\begin{figure}
\begin{center}
\begin{tabular}{l r}
  \includegraphics[scale=0.3]{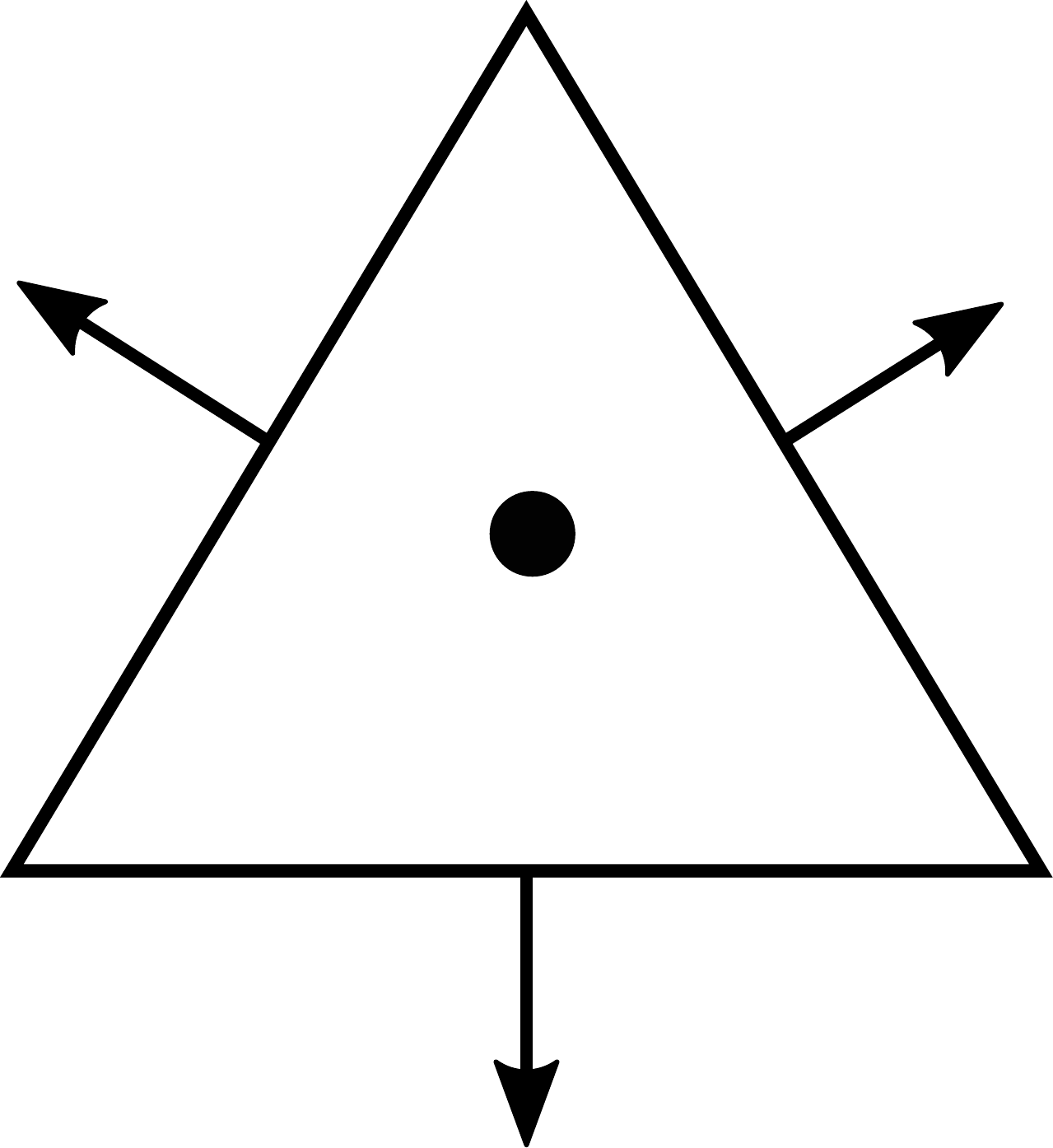}   & \includegraphics[scale=0.3]{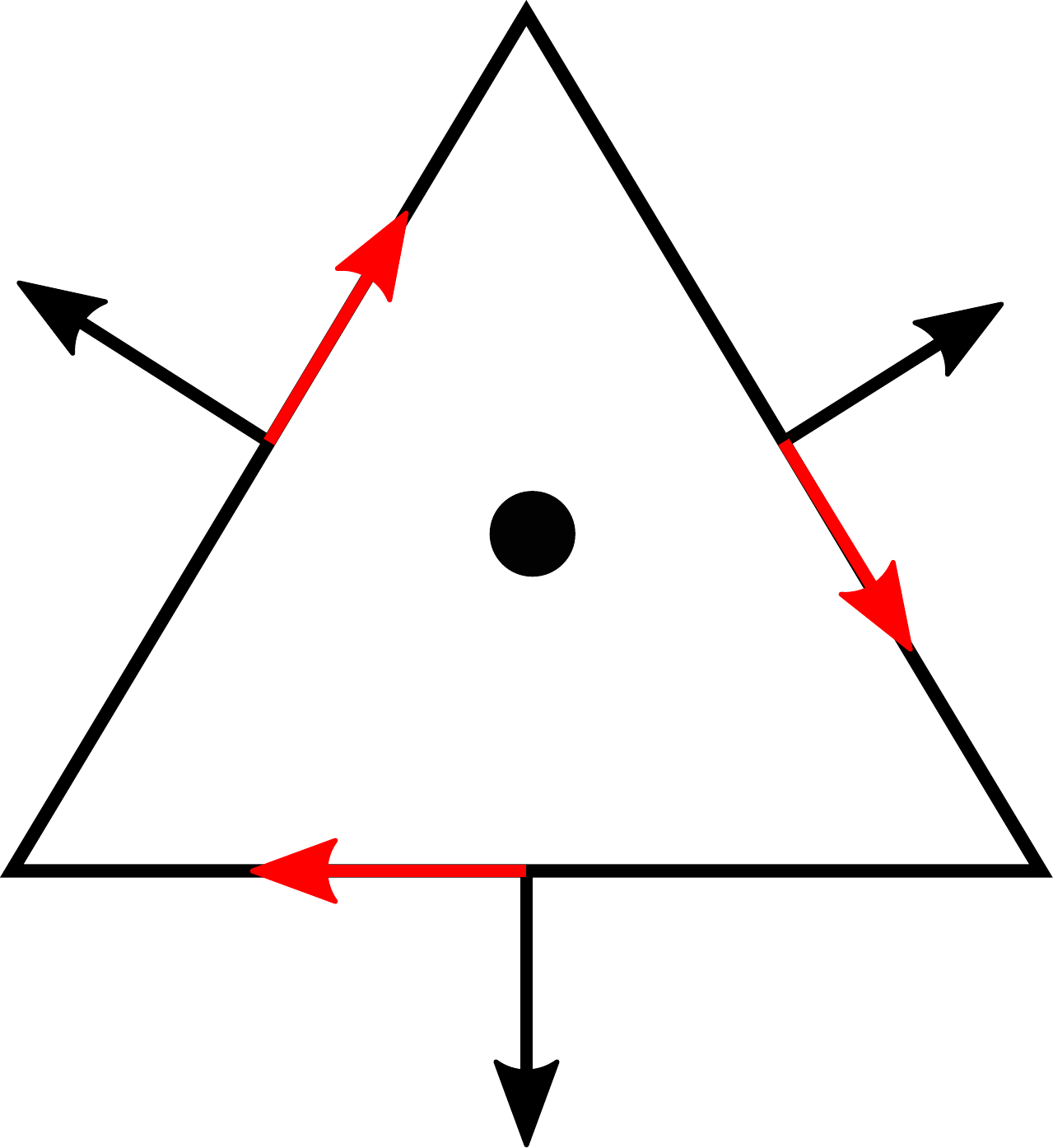} \\
\end{tabular}
\end{center}
\caption{{\it Left}: C-grid staggering of the ocean model ICON-O, where the normal {components} of horizontal velocities $\vt \cdot  \nt$  are placed at edge midpoints and the scalars are saved at cell centers. {\it Right}: The  Crouzeix-Raviart finite element,  where the horizontal velocity is represented by its normal and tangential component.\label{CCgrid}}
\end{figure}

\paragraph{Spatial discretization}
In the following we describe the { implementation} of the stabilized Crouzeix-Raviart finite element in ICON-O. 
Figure  \ref{CCgrid} shows the Raviart-Thomas element that corresponds to the C-type staggering in ICON-O and
the Crouzeix-Raviart finite element (CR). Both elements are coupled to a piecewise constant element for the scalar variables.

\begin{figure}
  \begin{center}
  \begin{tabular}{l c r}
      \includegraphics[scale=0.6]{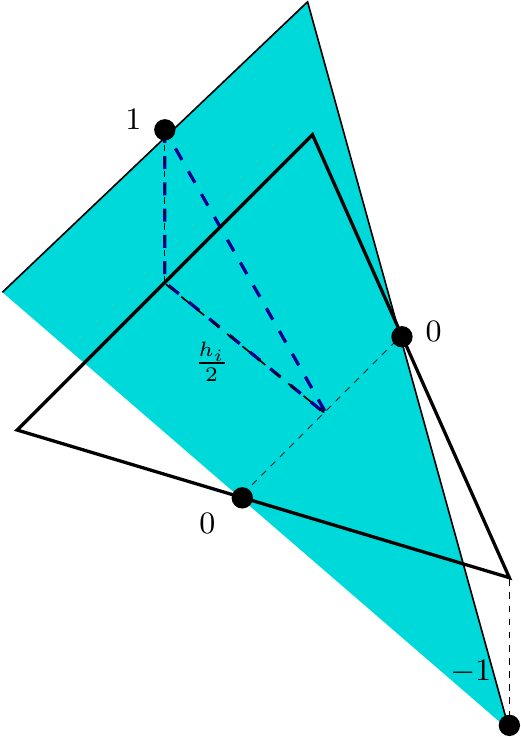}& &
    \includegraphics[scale=0.25]{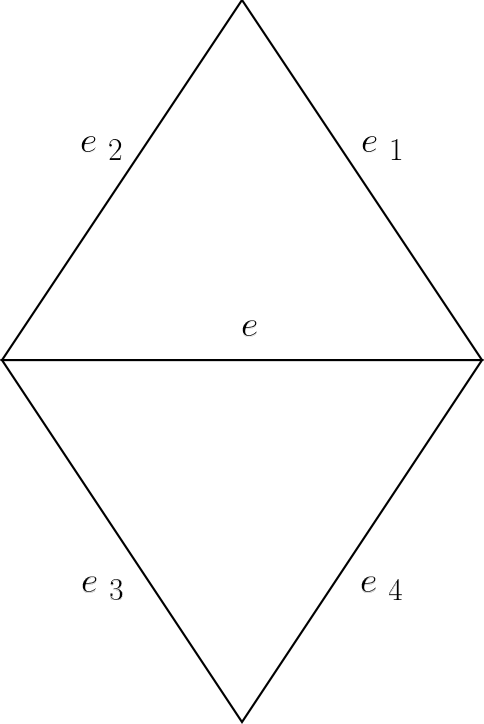}
  \end{tabular}
  \caption{{\it Left:} Basis function $\phi_i$  of the Crouzeix-Raviart finite element with $\phi_i(E_j)=\delta_{ij}$, where $E_j(K)$ refers to the edge midpoints of triangle $K$. {\it Right:} {\GG The support of a basis function at edge $e$ is given by {the two adjacent triangles. }}\label{testfunction}}
      \end{center}
  \end{figure}

By $\Omega_h$ we introduce the triangulation of a domain $\Omega$ into triangles $K$ that satisfy the usual assumption of structure and shape regularity. Let $V_h=\text{span} \, \{\phit_i, i=1,...,N\} \subset L^2(\Omega)$ be the space of the Crouzeix-Raviart element.  The edge mid points are denoted by $E_i\in \Omega_h$ and the edge itself by $e_i$. In each node $E_i$ a finite element basic function is given as  
\begin{align*}
\phi_i(E_j)&=\delta_{ij}, \quad  \forall i,j=1, ..,N, \quad \partial_n \phi_i(E_i)=\frac{2}{h_i},\quad  \partial_t \phi_i(E_i)= 0,
\end{align*}
where $N$ represents the number of edges in $\Omega_h$ and $h_i$ is the height orthogonal to edge $e_i$ of a triangle. To simplify the notation we define $\phi_i=\phi_i(E_i)$. The outward normal vector {\RR of the triangle} and the tangential vector to an edge $e_i$ at node $E_i$ is denoted by $\nt_i$ and $\taut_i$. The transposed is indicated by {\B the subscript} $T$. 
In Figure \ref{testfunction} we motivate that $\partial_{\tau} \phi_i= 0$ as the solid  line through 1 is constant.  Further, the differential quotient along the dashed blue lines in Figure  \ref{testfunction} gives  $ \partial_n \phi_i=\frac{1}{2^{-1}h_i}=\frac{2}{h_i}$. Using the discrete space $V_h$ we can express the velocity vector $\vt_h$ as
\begin{equation}\label{notation}
\begin{aligned}
\vt_h&=\sum^N_{i=1} (v_i\nt_i+u_i\taut_i) \phi_i,
\end{aligned}
\end{equation}
where $v_i$ and $u_i$ are scalar coefficients. 
The gradient of $\vt_h$ is given as
\begin{align*}
\nabla \vt_{h}=\sum^N_{i=1}(v_i\nt_i+u_i\taut_i) \nabla^T \phi_i=\sum^N_{i=1} (v_i\nt_i+u_i\taut_i)\nt_i^T\partial_n \phi_i,
\end{align*}
and the transposed gradient $\nabla \vt_h^T$ reads
\begin{align*}
\nabla \vt^T_{h}=\sum^N_{i=1}(v_i\nt_i^T+u_i\taut_i^T) \nabla \phi_i=\sum^N_{i=1} (v_i\nt_i^T+u_i\taut_i^T)\nt_i\partial_n \phi_i.
\end{align*}
We also decompose the basis functions $\phit$  and $\nabla \phit$ into the normal and tangential component 
\begin{align*}
\phit= \sum^N_{j=1} \nt_j \phi_j^n+\taut_j \phi_j^{\tau}, \quad \nabla \phit= \sum^N_{j=1}\nt_j \nt_j^T \partial_n \phi_j^n+ \taut_j\nt_j^T\partial_n \phi^{\tau}_j.
\end{align*}
Using these decomposition of the domain $\Omega_h$ in triangles, we can formulate the discrete momentum equation (\ref{eq:system}) as
{\R
\begin{equation}\label{eq:mom:spa:disc}
\begin{aligned}
 \sum_K &\sum_{E_i(K),E_j(K)}\Big(\rho h \partial_t (v_i\nt_i+u_i\taut_i) \phi_i -F(v_i\nt_i+u_i\taut_i)\phi_i, \nt_j \phi_j^n+\taut_j \phi_j^{\tau} \Big)_K\\
&+\Big (\sigmat_{i,h},\nt_j \nt_j^T \partial_n \phi_j^n+ \taut_j\nt_j^T\partial_n \phi^{\tau}_j \Big)_K 
= 0,
\end{aligned}
\end{equation}
}
where the subscript $K$ denotes the $L^2$-integral over a triangle $K$ and $E_i(K),E_j(K)$ are the edge midpoints of triangle $K$.
All integrals except of the one with stress tensor $ \Big(\sigmat_{i,h}, \nt_j \nt_j^T \partial_n \phi_j^n+ \taut_j\nt_j^T\partial_n \phi^{\tau} \Big)_K$  
{\R include the diagonal mass matrix. Thus, we can directly evaluate the sum over the edge midpoints in (\ref{eq:mom:spa:disc}). On each cell $K$ the contribution to the mass matrix is given as}
\begin{align*}
 M^{i,j}_K=\begin{pmatrix}
(\nt_i\phi^n_i,\nt_j\phi_j^n)_K& (\taut_i\phi^\tau_i,\nt_j\phi_j^n)_K\\
(\nt_i\phi^n_i,\taut_j\phi_j^\tau)_K& (\taut_i\phi^\tau_i,\taut_j\phi_j^\tau)_K
\end{pmatrix}=
\begin{pmatrix}
\delta_{i,j} \frac{|K|}{6}& 0\\
0&\delta_{i,j} \frac{|K|}{6}
\end{pmatrix},
\end{align*}
where $\delta_{i,j}$ is the Kronecker symbol and {\B $|K|$ the area of the triangle $K$}.
{\R As each edge is shared by two triangles the global mass matrix $M$ is diagonal
\begin{align*}
 M_{k,l}=\begin{pmatrix}
\delta_{k,l} \frac{|K|}{3}& 0\\
0& \delta_{k,l} \frac{|K|}{3}
\end{pmatrix},
\end{align*}
where $k,l=1...,N$.}
Next we consider the integral over the discretized stress tensor $\sigmat_h=\begin{pmatrix}
\sigmat_{11,h}&\sigmat_{12,h}\\
\sigmat_{12,h}&\sigmat_{22,h}\\
\end{pmatrix}$, which can be written as  
\begin{align*}
(\sigmat_h, \nabla \phit)_K=&\sum_{E_i(K), E_j(K)}(\sigmat_{i,h} \nt_j \nt_j^T \partial_n \phi_j^n+ \taut_j\nt_j^T\partial_n \phi_j^{\tau})_K\\
=&\sum_{E_i(K), E_j(K)}(\sigmat_{i,h} \nt_j \nt_j^T \frac{2}{|h_j|}+ \taut_j\nt_j^T\frac{2}{|h_j|})_K.
\end{align*}
 We evaluate the integral 
 and get  
\begin{align*}
 (\sigmat_{i,h},\nt_j \nt_j^T \partial_n \phi_j^n)_K = \sum_{E_i(K),E_j(K)}|K| \frac{2}{|h_j|}(\sigmat^{11}_{i,h} \nt^1_j\nt^1_j+2\sigmat^{12}_{i,h} \nt^1_j\nt^2_j+\sigmat^{22}_{i,h} \nt^2_j\nt^2_j),
\end{align*}
where we calculate the entries $\sigmat_{i,h}$ based on the strain rate tensor 
\begin{align*}
\frac{1}{2}(\nabla \vt_h + \nabla \vt_h^T)_T = \sum_{E_i(K)}(v_i\nt_i+u_i\taut_i)\nt_i^T \frac{2}{|h_i|}+(v_i\nt_i^T+u_i\taut_i^T)\nt_i \frac{2}{|h_i|}.
\end{align*}
The transport equations are discretized with an upwind scheme. 

\paragraph{Stabilization}

{\BR It is well-documented in the literature
that a discretization of the strain rate tensor with the Crouzeix-Raviart element causes a non-trivial null space in addition to the trivial null space generated by solid-body rotations (see  \cite{Falk1991} and \cite{Knobloch2000,Knobloch2004}). 
For an elasticity problem whose solution is unique up to the addition of the three-dimensional space of rigid motion, Falk \cite{Falk1991} showed that the corresponding discrete problem has a too large solution space if Korn's second inequality (\ref{Korn:second}) is not satisfied by the nonconforming piecewise linear finite element.
}
{\RR 
In case of pure Dirichlet boundary conditions the trivial null space is filtered out. A discretization of the strain rate tensor with the CR element can still cause a nontrivial kernel. Knobloch showed that for the Crouzeix-Raviart element Korn's first inequality (\ref{Korn:first}) is not uniform in the mesh size in the presence of Dirichlet boundary conditions.
To be precise,  {\GG  if $\vt_h$ denotes the discrete velocity of the Crouzeix-Raviart element, then it holds 
$\|\nabla \vt_h +\nabla \vt_h^T \|\ge c_k(h)\|\nabla \vt_h\|$ where the lower bound $c_k(h)$ depends on the resolution
such that $c_k(h)\to 0$ for $h\to 0$, see~\cite{Knobloch2004}. This is not sufficient to give robust approximations on fine meshes.   
}}
{\RR Brenner \cite{Brenner2004} derived a \emph{generalized} version of Korn's inequality
\begin{align}\label{eq:genkorn}
\| \nabla \vt_h \|^2 \leq c \Big (\| \nabla \vt_h+\nabla \vt_h^T \|^2 + \sum_{e} \int_{e} \frac{1}{|e|}[\vt_h][\vt_h] \, ds\Big), 
\end{align} 
where $[\vt_h]$ denotes jump of the velocity at an edge $e$ {\BR and is} defined by 
\begin{align}\label{jump:DG}
[\vt_h]:=\vt_h^+-\vt_h^-, \quad 
\vt_h^{\pm}:=\lim_{\epsilon\to 0} \vt_h^{\pm}(x \pm \epsilon \nt).
\end{align}
The generalized Korn inequality (\ref{eq:genkorn}) is fulfilled by the CR element. }
{\RR Motivated by the generalized inequality, we add to the momentum equation (\ref{eq:mom:spa:disc}) at each edge $e$ the stabilization  
\begin{align}\label{stabi_VP}
S_{i,j}^e:=2\zeta_{e} \frac{\alpha}{|e|} \int_{e} [\vt_{e_i}][\phit_{e_j}] \, ds,
\end{align}
with $i,j=1,..N$. {\BR In our work we  chose $\alpha=1$}. The stabilizing term can be interpreted as a discrete Laplacian, which penalizes the discontinuities along an edge. A similar stabilization was derived by Hansbo and Larson for a linear elasticity problem \cite{HansboLarson2003}. 
{\BB Adding the stabilization to the variational formulation  allows us to apply the generalized Korn inequality. Thus, the weighted 
$H^1$-estimates derived in Section \ref{sec:energy} are fulfilled by the stabilized CR approximation and the quantity $E(\vt)$ given in (\ref{eq:energy}) is bounded. 
{\BR This is a requirement for a well posed problem in the sense of Hadamard.} The numerical results in Figure \ref{fig:VPCR} show that without stabilization the approximation is not bounded and grows for $h \to 0.$}}

 
 Adding the stabilization to spatially discretized equation (\ref{eq:mom:spa:disc}) {\R yields the following modified momentum equation
\begin{equation}\label{eq:mom:stabi}
\begin{aligned}
&\sum_K \sum_{E_i(K), E_j(K)} \Big(\rho h \partial_t (v_i\nt_i+u_i\taut_i) \phi_i -F(v_i\nt_i+u_i\taut_i)\phi_i, \nt_j \phi_j^n+\taut_j \phi_j^{\tau} \Big)_K\\
&+\Big (\sigmat_{i,h},\nt_j \nt_j^T \partial_n \phi_j^n+ \taut_j\nt_j^T\partial_n \phi_j^{\tau} \Big)_K+ \sum_e \sum^N_{i,j=1}  S^e_{i,j}
= 0.
\end{aligned}
\end{equation}}
{\R For the calculation of the stabilization term in (\ref{eq:mom:stabi}) we have to evaluate the integrals {\R $\int_{e} [\vt_{e_i}][\phit_{e_j}] $} in (\ref{stabi_VP}).
For computing these integrals 
 \begin{align}\label{eq:jumps}
 \int_{e} [\vt_{e_i}][\phit_{e_j}] ds
 \end{align}} one has to take into account the coupling of the  test and ansatz functions along the five neighboring edges shown in Figure \ref{testfunction}.  Since $\int_{e} [\phit_e]\, ds=0$, the stencil reduces to the test and ansatz functions defined at the four surrounding  edges $e_{j}$, $j=1,..,4$.  For $i,j=1,...,4$ the integral over an edge $e$ is given as 
\begin{align*}
\int_{e}\phi_{i}\phi_{j}\, ds=\begin{cases} &\frac{1}{3}|e|,\,\text{if }\,  i=j,\\
 -&\frac{1}{3}|e|, \, \text{else},
 \end{cases}
\end{align*}
where $|e|$ is the length of edge $e$.
 Next, we reformulate expression (\ref{eq:jumps}) and get for $i,j=1,...,4$

\begin{align*}
\int_{e} [\vt_{e_i}][\phit_{e_j}] \,  ds
= \sum_{E_{i}, E_{j}} \int_{e}[v_{i} \nt_{i} \phi_{i}^n+u_{i}\taut_{i} \phi_{i}^\tau][\nt_{j}\phi_{j}^n+\taut_{j} \phi_{j}^\tau] ds,
\end{align*} 
where $v_{i}$ and $u_{i}$ are the coefficients of the velocity located at the midpoint $E_{i}$ of an edge.
{\B The integral  $\int_{e} [\vt_{e_i}][\phit_{e_j}] $ can be efficiently evaluated, by two loops over all triangles. }

\section{Numerical evaluation}\label{sec:num}
This section provides an experimental analysis of the discretization of the sea-ice momentum equation with the Crouzeix-Raviart element. Section \ref{sec:strain} starts with analyzing the strain rate tensor in the viscous regime. In Section \ref{sec:VPEVP} we study the full viscous-plastic and elastic-viscous-plastic rheology and investigate {\RR if the CR approximation fulfills the $H^1$-estimate defined in (\ref{eq:est})}. Finally in Section \ref{sec:box}, we evaluate the full system describing the sea-ice dynamics including the advection of the mean sea-ice thickness and the sea-ice concentration. We analyze a box test, which is a slight modified version of the test case introduced by Danilov et al.~\cite{Danilov2015}. 
\subsection{Strain rate tensor}\label{sec:strain}
\begin{figure}
\begin{tabular}{c c c c }
 \includegraphics[scale=0.7]{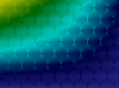}&  \includegraphics[scale=0.7]{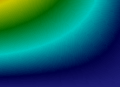}&  \includegraphics[scale=0.7]{pics_CR/epsilon_mS_z1.png}\\
  \includegraphics[scale=0.2]{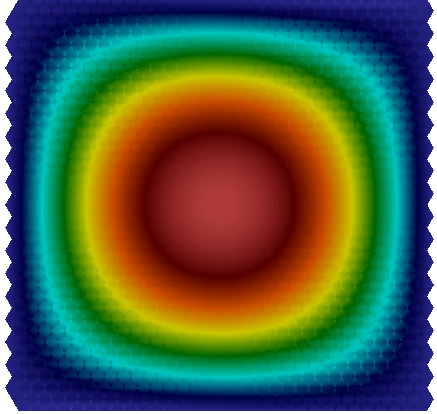}&  \includegraphics[scale=0.2]{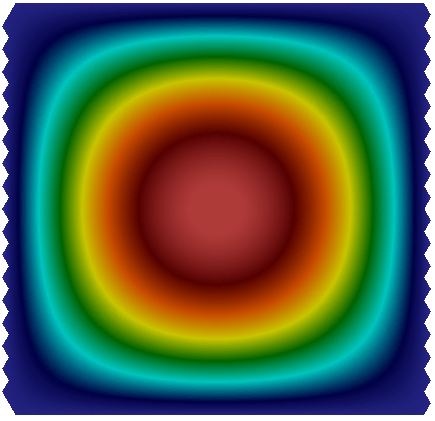}&  \includegraphics[scale=0.2]{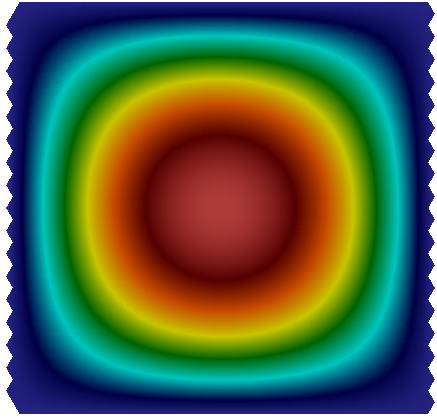}&
\includegraphics[scale=0.2]{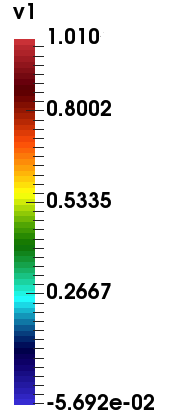}\\
$\sigmat=\frac{\zeta}{2}\Big( \nabla \vt +\nabla \vt^T \Big )$  & $\sigmat_1=\frac{\zeta}{2}\nabla \vt$ & $\sigmat$ stabilized&\\
\end{tabular}
\caption{We consider the simplified form of the momentum equation (\ref{eq:mo:sim}) with different stress tensors. 
The lower row shows the first velocity component and the upper row presents a close-up from the lower right corner. The left column
shows the solution for the Crouzeix-Raviart element. The middle figures refer to the case of the reduced stress tensor 
$\sigmat_1$. The plot on the right shows the solution for the full stress tensor $\sigmat$ with the stabilized Crouzeix-Raviart element.  \label{fig:strain}}
\end{figure}

 We start  our analysis with a simplified version of the momentum equation (\ref{eq:mo}).
\begin{align}\label{eq:mo:sim}
\partial_t \vt -\div( \sigmat(\vt))= \rhs,
\end{align}
and consider the viscous-plastic stress tensor in the viscous regime. 
 \begin{align*}
 \sigmat = \zeta\frac{1}{2}(\nabla \vt +\vt^T),\quad \zeta=\frac{P}{2\Delta_\text{min}}, \quad h=1,\, A=1.
\end{align*}
The advection of the mean sea-ice thickness and sea-ice concentration is deactivated and we set $h=1,\, A=1$. The domain is a planar quadrilateral with length of $L_x=L_y=500$ km in $x$ and $y$ direction and tessellated by a mesh of equilateral triangles. We  start the simulation with zero initial velocities and apply homogeneous Dirichlet conditions at the boundaries. 
Given an analytic solution $\vt^1=\vt^2=-\sin\Big({\pi_x x}\Big)\sin\Big({\pi_y y}\Big)$, with $\pi_x:=\frac{\pi}{L_x}$ and   $\pi_y:=\frac{\pi}{L_y}$, the right hand side 
of (\ref{eq:mo:sim}) is
\begin{align*}
\rhs=\frac{\zeta}{2} \Big({\pi_x^2}\sin({\pi_x} x)\sin({\pi_y} y) &+\frac{\pi^2_y}{ 2}\sin({\pi_x} x)\sin({\pi_y}y)\\
&-\frac{1}{2}{\pi_x}{\pi_y}\cos({\pi_x} x)\cos({\pi_y} y)\Big).
\end{align*}
We observe instabilities in the velocity field shown in left plot in Figure \ref{fig:strain}. 
This behaviour is consistent with fact that the Crouzeix-Raviart {\RR element fulfills} {\G the discrete version of Korn's inequality (\ref{Korn:first})
 \begin{align*}
 c_k(h) \| \nabla \vt_h\|^2\leq \| \nabla \vt_h +\nabla \vt_h^T \|,
 \end{align*}
only with a positive mesh depended constant $c_k(h)$.
}
This instabilities vanish if we consider $\sigmat_1=\frac{\zeta}{2}\nabla \vt$  instead of the symmetric stress tensor $\sigmat=\frac{\zeta}{2}\Big ( \nabla \vt +\nabla \vt ^T \big )$.
The oscillations also disappear if we add the stabilization {\R $S^e_{i,j}$ (\ref{stabi_VP}) to the momentum equation (\ref{eq:mo:sim}).}
We observe that the stabilization slightly damps the solution. 
Without stabilization the velocity components reach their maximum at  $\vt_h^1=1.052$ and $\vt^2_h=1.28$. Replacing $\sigmat$ by $\sigmat_1$ reduces the maximal velocity to $\vt_h^1=\vt^2_h=1.027$. Adding the edge-stabilization for $\sigmat$ to equation (\ref{eq:mo:sim}) the velocities further  decrease the maxima  to $\vt^1_h=1.023$ and $\vt^2_h=1.005$. For the runs that include $\sigmat_1$ we adjusted the right hand side of equation (\ref{eq:mo:sim}) to $\rhs=\zeta\Big ({\pi^2_x}+{\pi^2_y}\Big)\sin\Big(\pi_x x \Big)\sin\Big(\pi_y y\Big)$ to converge against the same analytic solution $\vt^1=\vt^2=-\sin\Big({\pi_x x}\Big)\sin\Big({\pi_y y}\Big)$.  All simulations presented in Figure \ref{fig:strain} are computed with an explicit Euler method using a time step $k=1 \cdot 10^{-6}$ s  on a triangular mesh with 3833 edges. 
\subsection{Viscous-plastic and elastic-viscous-plastic rheology}\label{sec:VPEVP}
We numerically analyze the effect of the Crouzeix-Raviart discretization on the approximations of the EVP and VP model.  For our investigation we neglect the wind stress and the Coriolis force and consider 
\begin{align}\label{sec:VPEVP_test}
\rho h\partial_t \vt=\div(\sigmat)-\rho_wC_{w}\|\vt-\vt_{w}\|(\vt-\vt_{w}).
\end{align}
Following Hunke \cite{Hunke2001} and Danilov et al.~\cite{Danilov2015} we defined the ocean velocity as
\begin{align*}
\vt^1_w=\frac{0.1(2y-L_y)}{L_y},\quad \vt^2_w= \frac{0.1(L_x-2x)}{L_x}.
\end{align*}
The viscous-plastic stress tensor and the elastic-viscous-plastic stress tensor is given in (\ref{stress}) and in (\ref{EVP}). 
We consider the same quadrilateral domain as in Section \ref{sec:strain} with a triangular grid and assume homogeneous Dirichlet boundary conditions. The initial ice velocity is $\vt(t_0)=0$.

In a first test we switch off the  advection of the mean sea-ice thickness and sea-ice concentration and use $h=1$ and $A=\frac{x}{L_x}$. For this configuration the velocities convergence against a stationary solution, such that we can compare the solutions of the EVP model to those of the VP model. The approximation of the momentum equation with the viscous-plastic rheology is solved  in two ways. First, fully explicit, using the forward Euler time-stepping method and a time step of $k=0.1$ s, second, with the mEVP solver described in Section \ref{sec:disc} and a time step of $k=600$ s and $N_\text{mevp}=500$ sub-iterations. For both cases we observe instabilities in the velocity field. The same holds for the EVP model. 
\begin{center}
\begin{figure}
\begin{tabular}{c c c  }
 \includegraphics[scale=0.15]{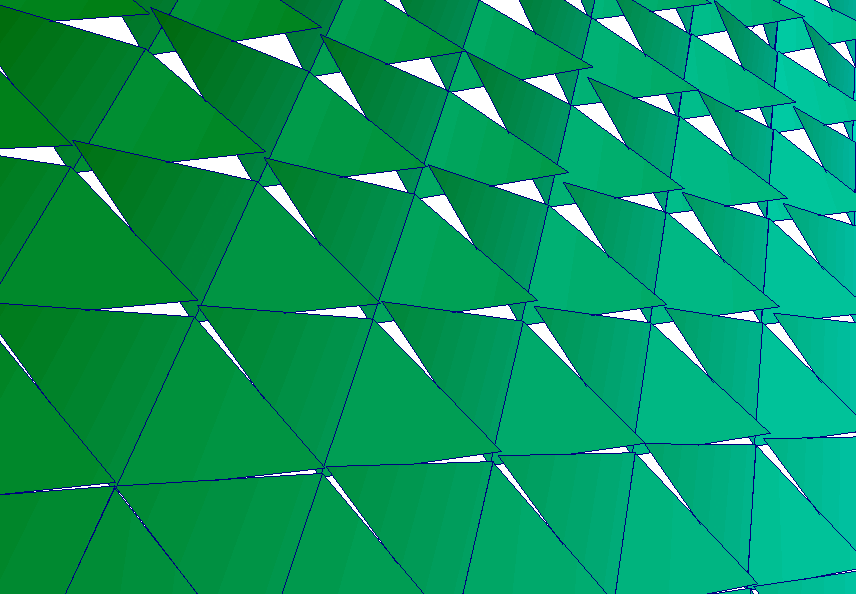}&  \includegraphics[scale=0.15]{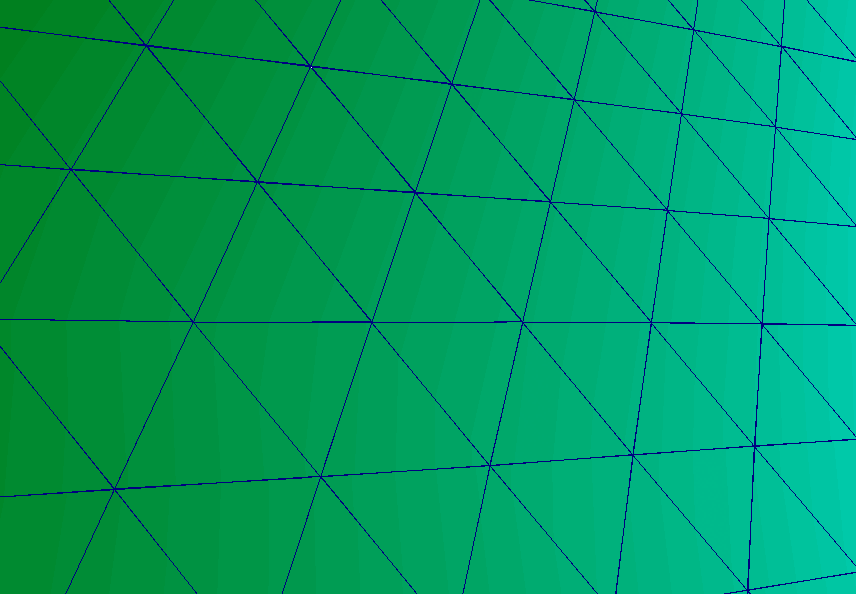}&\\
  \includegraphics[scale=0.3]{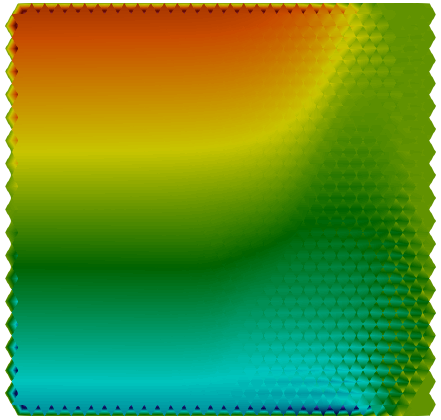}& 
   \includegraphics[scale=0.3]{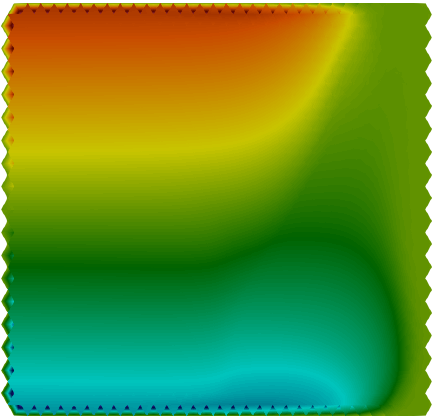}&
   \includegraphics[scale=0.3]{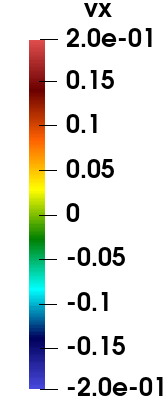}\\
   VP &stabilized VP &\\
  \includegraphics[scale=0.3]{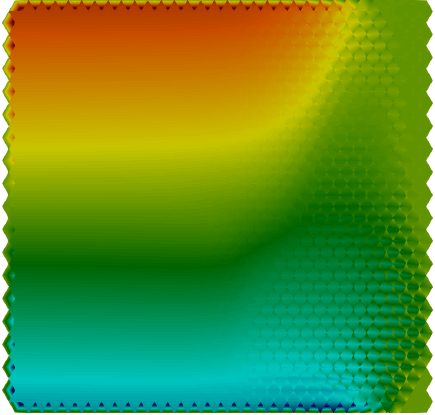}& 
  \includegraphics[scale=0.3]{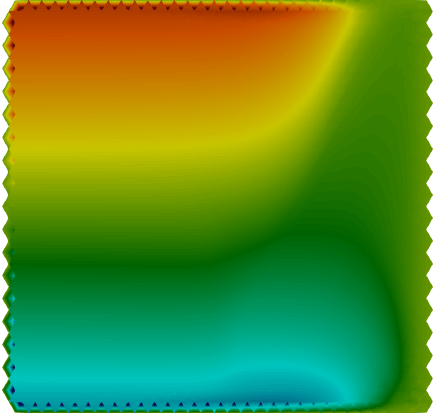}&
   \includegraphics[scale=0.3]{pics_CR/review_stat_scalarbar.png}\\
   mEVP &stabilized mEVP &\\
  \includegraphics[scale=0.3]{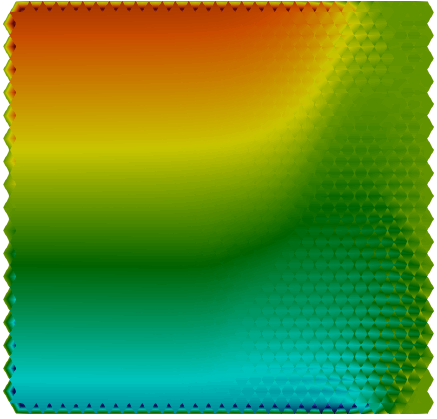}&
   \includegraphics[scale=0.3]{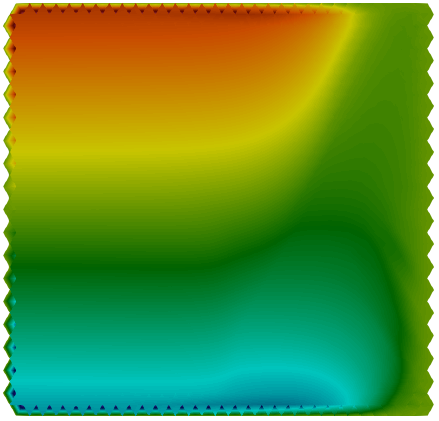}&
 \includegraphics[scale=0.3]{pics_CR/review_stat_scalarbar.png}\\
EVP &stabilized EVP &\\
\end{tabular}
\caption{\small{The first component { $\vt^x$ of the velocity vector} of test case (\ref{sec:VPEVP_test}). We evaluate the VP, EVP { model} and mEVP { approximation} without stabilization (left column) and with stabilization (right column) of the momentum equation. In the bottom row we zoom into the lower right corner of  the velocity calculated with the mEVP solver.\label{fig:vp_evp}}}
\end{figure}
\end{center} 
As shown in Figure \ref{fig:vp_evp} instabilities appear in the velocity field in regions with high sea-ice concentration. With the corresponding stabilization these instabilities vanish and all three approximations produce similar results.
We zoom into the approximation computed with the mEVP solver and show the lower right corner of the domain. One can see that using the edge-stabilization the oscillation in the Crouzeix-Raviart element completely disappears. 

{\G The {\BR weighted $H^1$-norm} 
\begin{align}\label{Weighted_Hone}
E(\vt)= \int_I \zeta_\text{min}\|\nabla \vt\|^2 dt 
\end{align}
was defined in (\ref{eq:energy}) in Section \ref{sec:energy}. 
{\BR  
Figure \ref{fig:VPCR} shows $E(\vt)$ {\BR for the stabilized as well as for the}  unstabilized CR approximation. The upper panel in Figure \ref{fig:VPCR} shows $E(\vt)$ on a 10 km mesh as a function of time. We evaluate $E(\vt)$ for the EVP model and the VP model with and without stabilization. In case of the VP model we use either a fully explicit discretization in time or the mEVP solver.  The different configurations presented in Figure \ref{fig:VPCR} correspond to the different discrete solutions presented in Figure \ref{fig:vp_evp}. From the upper panel of Figure \ref{fig:VPCR} one can infer that
$E(\vt)$ increases linear with time as expected, and that the unstabilized approximations grow faster than the stabilized one. {\GG The behavior of $E(\vt)$ is relatively similar within the stabilized approximations and within the unstabilized configuration.} Therefore we restrict our further analysis of the difference between stabilized and non-stabilized case to the representative mEVP model approximation.


The middle panel in Figure \ref{fig:VPCR} shows in the mEVP case $E(\vt)$ for three different resolutions: 10 km, 5 km and 2.5 km.
For the stabilized case no change in $E(\vt)$ is visible under a resolution increase, while the unstabilzed  CR approximation
shows not only differences for $E(\vt)$ at different resolutions but in fact a higher growth rate with higher resolution.  

This is further illustrated in the lower panel of Figure \ref{fig:VPCR},  which shows $E(\vt)$ in a logarithmic scale on a 10 km, 5 km and 2.5 km mesh at a fixed point in time $T=24$ h. We observe that  in the stabilized configuration $E(\vt)$ stays almost constant, whereas
in the unstabilized case $E(\vt)$ grows linear under resolution increase.


The weighted $H^1$-estimate in (\ref{eq:energy}) for {\GG the VP and EVP rheologies state that $E(\vt)$ of the solution of the continuous sea-ice equations is bounded by the data. This property should be shared by the discrete  approximations at arbitrary mesh resolutions.}
Let us denote the upper bound of the estimate {\GG at a fixed instant in time} by $R$. The lower panel in Figure \ref{fig:VPCR} suggests that for decreasing mesh size $h \to 0$ 
the quantity $E(\vt)$ exceeds the upper bound $R$.} {\GG Thus, we conclude that the Crouzeix-Raviart element without stabilization has a qualitative different behaviour than the solution of the continuous equations.}
{\BB 
This observation is consistent with the result from Knobloch \cite{Knobloch2004}, where it is shown that a mesh dependent constant $c_k(h)$ prevented the Crouzeix-Raviart element from the convergence against the solution of a Stokes problem coupled to deformation including the strain rate tensor $\epsilont$.}

 \begin{figure}
 \includegraphics[scale=0.65]{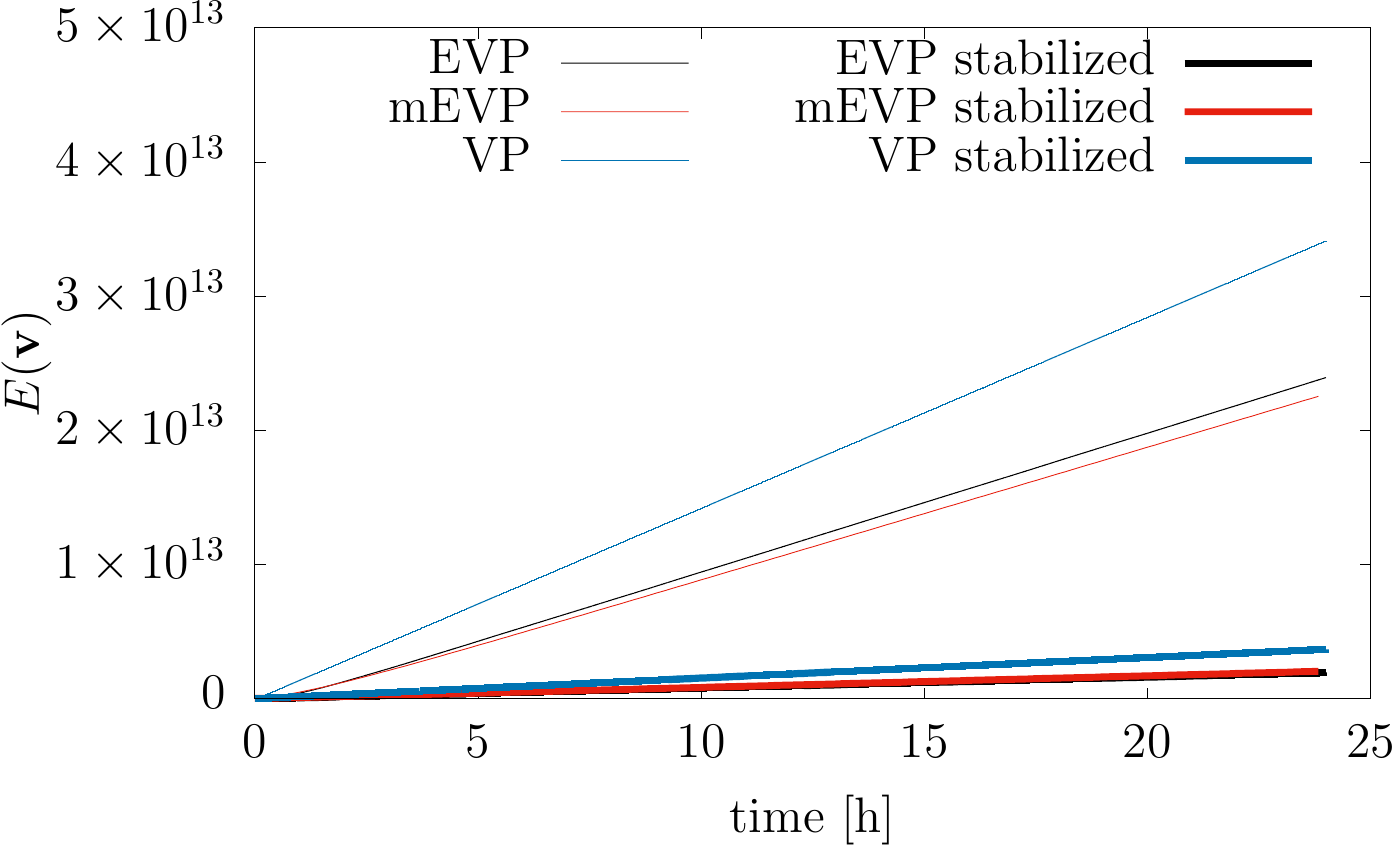}
 \includegraphics[scale=0.65]{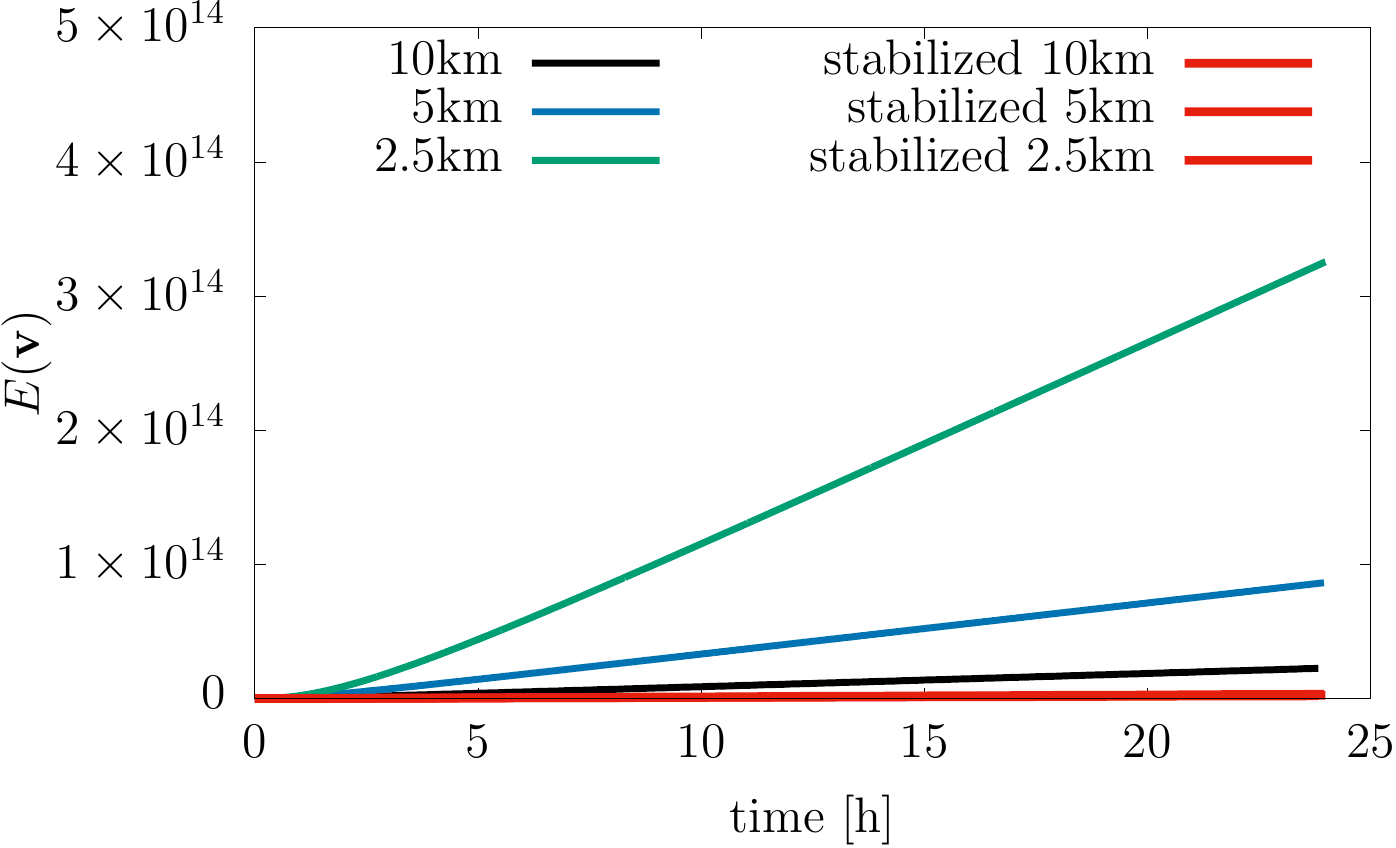}
 \includegraphics[scale=0.65]{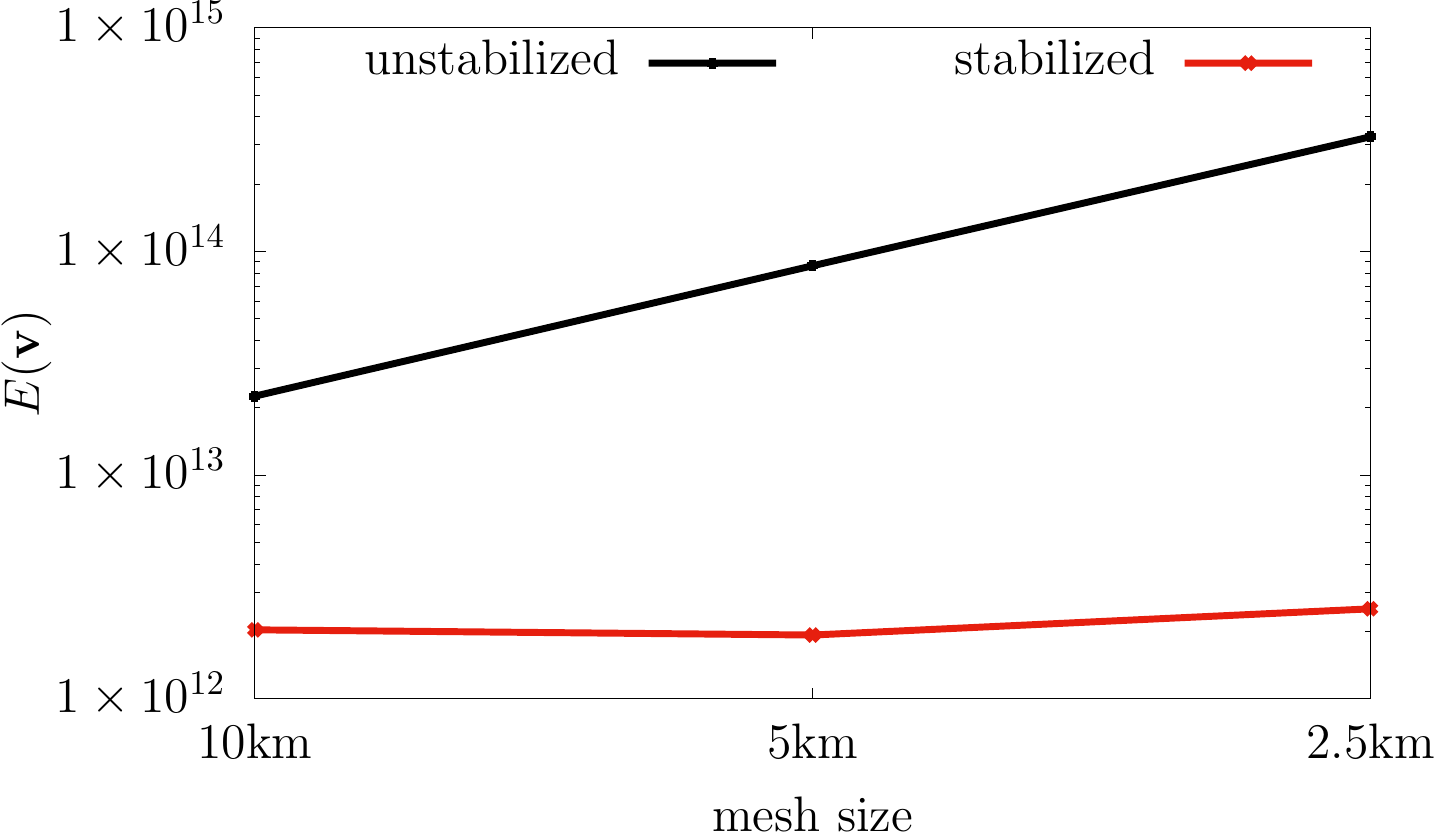}
  {\small\caption{ $E(\vt)$, defined in (\ref{Weighted_Hone}), computed with stabilized and unstabilized Crouzeix-Raviart element. 
   {\it Top:} $E(\vt)$ as a function of time for EVP, VP and mEVP. {The three unstabilized solutions show a qualitatively similar linear increase of $E(\vt)$ with time. The same applies to the three stabilized solutions, but with a smaller growth rate.} {\it Middle:} $E(\vt)$ of mEVP {approximation} at different resolutions. The plots of the stabilized cases are indistinguishable, while the unstabilized ones show an increasing growth with increasing resolution. {\it Bottom}: $E(\vt)$ plotted as function of resolution at fixed time $T=24$ h.   \label{fig:VPCR}}}
 \end{figure}
 
 {\RR In a second test (not shown) with activated advection we observed that the instabilities of the velocities propagated into the tracers fields. }
 }
 

\subsection{Box test}\label{sec:box} 
In this section we investigate the full system of sea-ice equations with a test case that is a slightly modified version of the box test described by Danilov et al.~\cite{Danilov2015}. The domain is a square of length $L_x=L_y=1000$ km. It is discretized with a triangular mesh
of equilateral triangles with a side length of approximately 15 km and 15190 edges. The  ocean current is as described in Section \ref{sec:VPEVP}. The wind velocity is given  by
\begin{align*}
\vt^1_\text{a}=5+(\sin(2\pi t/T)-3)(\sin(2\pi x/L_x)\sin(2 \pi y/L_y),\\
\vt^2_\text{a}=5+(\sin(2\pi t/T)-3)(\sin(2\pi y/L_y)\sin(2 \pi x/L_x).\\
\end{align*} 
The considered time span is one month.
We assume homogeneous Dirichlet boundary conditions and initial data given by $\vt(t_0)=0$, $h(t_0)=1$ and $A(t_0)=\frac{x}{L_x}$. {A time step of {$k=600$ s and $N_\text{evp}=500$ sub-cycles} are used.} 
In a first experiment the momentum equation is solved with the mEVP solver with constant  $A=\frac{x}{L_x}$ and $h=1.0$.
\begin{figure}
\begin{center}
\begin{tabular}{c c c}
 \includegraphics[scale=1.3]{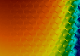}
  &  \includegraphics[scale=1.30]{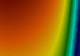}&  \\
  \includegraphics[scale=0.25]{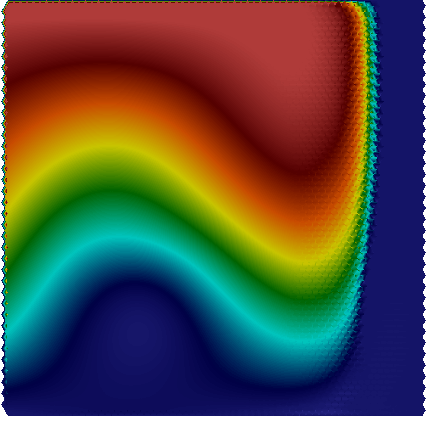}
  &  \includegraphics[scale=0.25]{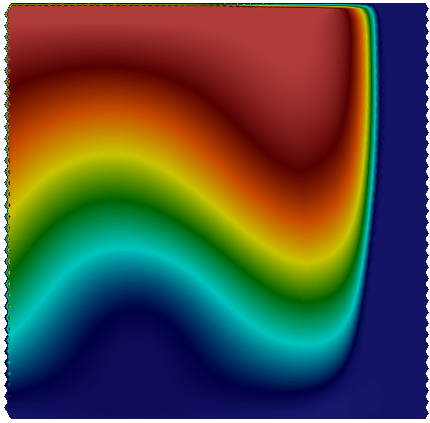}&  
     \includegraphics[scale=0.25]{pics_CR/review_stat_scalarbar.png}\\
      \includegraphics[scale=0.25]{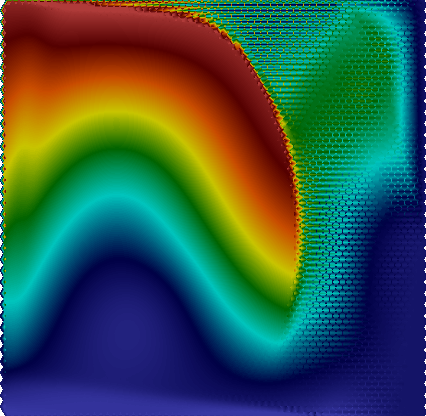}& \includegraphics[scale=0.25]{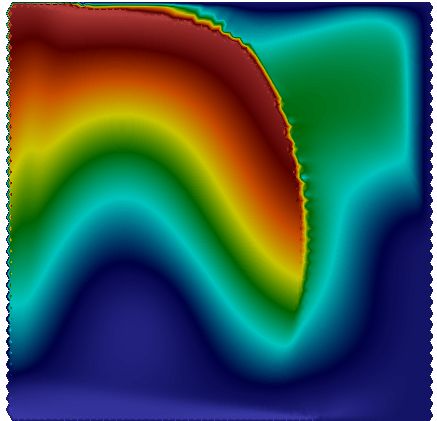}& \includegraphics[scale=0.25]{pics_CR/review_stat_scalarbar.png}\\
        \includegraphics[scale=0.68]{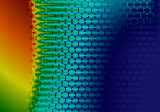}& \includegraphics[scale=0.68]{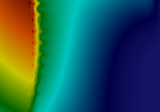}&\\
mEVP& stabilized mEVP& \\
\end{tabular}
 \caption{ {The first component $\vt^x$  of the sea-ice velocity vector without stabilization (left column) and with stabilization (right column). In the top row we show  $\vt^x$ after one month of simulation  with fixed sea-ice concentration and mean sea-ice thickness. The second row presents a close-up from the upper right corner. The third row refers to  $\vt^x$ after one week of simulation  with active advection. A closed up of the upper right corner is shown in the bottom row.\label{fig:boxteststat}}}
  \end{center}
\end{figure}

\begin{figure}
\begin{center}
\begin{tabular}{c c c}
\includegraphics[scale=1.2]{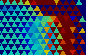}&
  \includegraphics[scale=1.2]{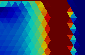}&\\
      \includegraphics[scale=0.25]{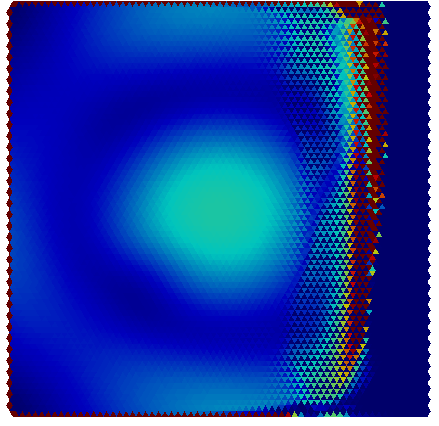}&
  \includegraphics[scale=0.25]{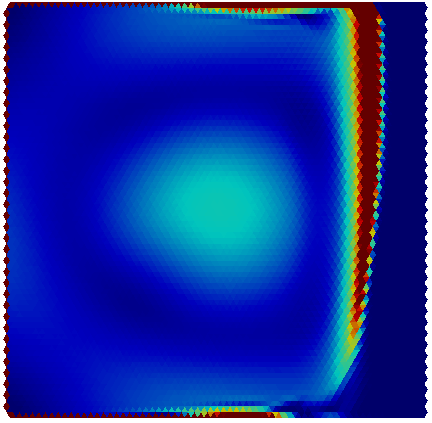}&
     \includegraphics[scale=0.26]{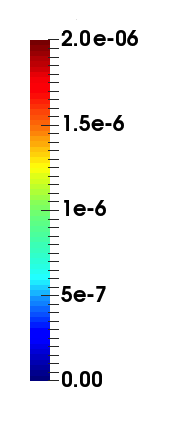}\\
mEVP& stabilized mEVP \\
\end{tabular} 
 \caption{{ The lower row shows the approximation of the total deformation} $\Delta$ { (cf. (\ref{viscosities})) with and without stabilization after one month of simulation with fixed sea-ice concentration and  mean sea-ice thickness\label{fig:VP_CR_AV}. The upper row presents a close-up from the upper right corner.}}
  \end{center}
\end{figure}

After one month we observe instabilities in the velocities in regions with large gradients in $\Delta$ (defined in (\ref{viscosities})) and high sea-ice concentrations. As shown in Figure \ref{fig:boxteststat} and Figure \ref{fig:VP_CR_AV}, the oscillations vanish if the stabilization is applied. 
With active advection we note also instabilities in the velocity field in regions with high sea-ice concentration. As shown in Figure \ref{fig:boxteststat} and Figure \ref{fig:VP_CR_AV} stabilizing the momentum equation results in a smooth representation of the velocities. 
{\B We stopped the simulation after one week. The underlying uniform mesh is too coarse to resolve the sharp gradient in the tracer that develops in the upper left corner of the domain. Here artefacts along the sea-ice edge started to appear. As we observe the same behaviour with a piecewise linear finite element discretization we {\RR attribute} those artefacts to the resolution of the underlying mesh and not the choice of the finite element.}


\section{Conclusion}
In this paper we introduced a new sea-ice discretization on the triangular grid in ICON. The discretization is based on a stabilized nonconforming Crouzeix-Raviart finite element, which consists of normal and tangential velocity components staggered at the edge midpoints of a triangle. The velocities are coupled to cell wise constant representations of the tracers. This staggering allows straightforward coupling to C-grid ocean and atmosphere discretizations. We numerically showed that a direct discretization with the Crouzeix-Raviart element leads to an unstable approximation of the velocities, which stems from the discretization of the strain rate tensor. 
To {\RR overcome} this issue we introduced an edge-based stabilization. We demonstrated numerically that stabilizing the sea-ice velocity is necessary for both the viscous-plastic and elastic-viscous-plastic model. {\GG To show that the approximation with the stabilized CR element is consistent with the solution of the continuous sea-ice equations, we derived a $H^1$-estimate for the VP and EVP model. {\R The $H
^1$-estimate bounds the {\BB weighted gradient $E(\vt)$ of the sea-ice system.}} We numerically evaluated $E(\vt)$ for an approximation of the viscous-plastic and elastic-viscous-plastic model and found that with the stabilization of the Crouzeix-Raviart element $E(\vt)$ stays bounded as in the continuous case. Without stabilization the {\BB this quantity} 
{\RR grows} with increasing mesh resolution, {\RR which is a qualitative different behaviour compared to solution of the continuous sea-ice equations}}. This underlines the importance of stabilizing the  Crouzeix-Raviart element when discretizing the sea-ice momentum equation. 
{\B We integrated the new sea-ice discretization in the coupled sea-ice ocean model in ICON. Due to the C-grid type staggering the new discretization benefits from a straight forward coupling to the C-grid ocean discretization.} {\RR The analysis of the stabilization on nonuniform grids is subject of future work.}

\begin{appendix}
{\RR
\section{Appendix}\label{app}}
\paragraph{Functional spaces and Korn's inequality}
{\G  The Sobolev space $H^1(\Omega)^2$ consist of functions, defined on 
a two-dimensional domain $\Omega\subseteq \mathbb{R}^2$ with square-integrable derivatives. The space $H^1_0(\Omega)$ contains all function of $H^1(\Omega)$ with trace zero on the boundary
$\Gamma$. We prescribe homogeneous Dirichlet conditions $\vt=0$
on the whole boundary  $\Gamma:=\partial\Omega$.  In order to describe time-varying function $\vt(x,t)$
we introduce the time interval $[0,T]$, $T>0$, and the space $L^2([0,T], H^1_0(\Omega)^2)$ that consists of functions such that the norm 
\begin{align*}
   ||\vt||_{L^2([0,T], H^1_0(\Omega)^2)}^2:=\int_0^T||\vt(t)||_{H^1_0(\Omega)^2}\, dt 
\end{align*} is finite. Based on the definition of the functional spaces we formulate Korn's {\BR first and second inequality} \cite{Ciarlet2010}.  Korn's first inequality states that there exists a positive constant $c_k$ such that 
\begin{align}\label{Korn:first}
\|\epsilont \|^2 \geq  c_k \|\nabla \vt \|^2 , \forall \vt \in H^1_0(\Omega). 
\end{align}
To introduce Korn's second inequality we define the space
{\RR
\begin{align}\label{Korn:space}
   W:=\{ \vt\in H^1(\Omega)| \int_\Omega \vt \,dxdy =0,\quad  \int_\Omega \text{rot}(\vt)  \, dxdy =0\}.
\end{align}
}
Now, Korn's second inequality states that for all $\vt \in W$ exits a positive constant $c_k$ such that 
\begin{align}\label{Korn:second}
\|\epsilont \|^2 \geq  c_k \|\nabla \vt \|^2 , \forall \vt \in W(\Omega).     
\end{align}
}

\begin{theorem}[$H^1$-estimate for viscous-plastic sea-ice momentum equation]\label{thoe1}
\hfill
{\B
Let a time interval [0,T], with $T>0$, be given. Let sea-ice concentration $A$ and { mean} sea-ice thickness $h$ be constant in time and space and{ $A$, $h>0$}.  Suppose
$\vt \in L^2([0,T],H^1_0(\Omega))$ 
is a solution of the sea-ice momentum equation (\ref{eq:mo}) with VP rheology. 
Then $\vt$
satisfies the following $H^1$-estimate 
\small{
\begin{align}\label{eq:estimate}
\|\rho h \vt(T) \|^2
&+ \int_{0}^T  \{ 2\rho_wC_{w1} \|\vt \|^2 \\
&+    c_k \frac{\zeta_\text{min}}{4} \|\nabla \vt \|^2 \} dt \leq \int_{0}^T { \frac{c_p}{ c_k\zeta_\text{min}} } \|\rhs\|^2 dt 
+\|\rho h \vt(0) \|^2,
\end{align}}
where $c_k$ and $c_p$ are positive constants that depend on the domain.{\GG The minimal value of the viscosity is defined as $\zeta_\text{min}:=\min_{(\xt,t) \in \Omega\times I} \zeta(\xt,t)$} and  $\rhs:=\rho_a C_a\|\vt_a\| \vt_a+\rho_wC_{w1}\vt_w-\rho h f\te \times \vt_{w}$. 
}
\end{theorem}
{\G
\begin{proof}
We consider the momentum equation in the weak formulation  (\ref{eq:system}) with $\phit=\vt$. 
\begin{align}\label{eq:weakest}
( \rho h\partial_t \vt, \vt) -(F+\div(\sigmat), \vt)=0.
\end{align}
Then, we apply the chain rule to the time dependent integral and get 
\begin{align}\label{eq:temp}
(\rho h\partial_t \vt, \vt)= \frac{1}{2} \int_{\Omega} \rho h \partial_t |\vt |^2 \, dxdy =  \frac{1}{2} \partial_t  \|\rho h \vt \|^2.
\end{align}
We proceed with analyzing the external forces
\begin{align}\label{eq:exforce}
(F,\vt)=( \tau+\rho h f\te \times \vt_{w}-\rho h f_c \te \times \vt,\vt).
\end{align} 
The integral over the Coriolis force $(f_c \te \times \vt,\vt)$ vanishes as the product $(\vec a\times \vec
b)\cdot \vec c=-(\vec a\times \vec c)\cdot \vec b$ is anti-commutative.
{The wind and ocean drag simplifies to 
\begin{align*}
\Big( \tau, \vt \Big)
&=(\rho_a C_a\|\vt_a\| \vt_a+\rho_w C_{w1}(\vt_w-\vt),\vt)\\
&=(\rho_a C_a\|\vt_a\| \vt_a+\rho_w C_{w1}\vt_w,\vt)-\rho_wC_{w1} \|\vt \|^2.
\end{align*}
}
{We define the $\vt$-dependent part of the ocean drag $(\tau,\vt)$
as 
\begin{align}\label{eq:oceandrag}
    \tau(\vt)\vt:=- \rho_wC_{w1} \|\vt \|^2,
\end{align}}
and collect all $\vt$-independent terms of equation (\ref{eq:exforce})
\begin{align*}
\rhs=\rho_a C_a\|\vt_a\| \vt_a+\rho_wC_{w1}\vt_w+\rho h f\te \times \vt_{w}.
\end{align*}
{This implies for the right hand side of equation {\G (\ref{eq:weakest})}}
\begin{align}\label{eq:to_be_estimated}
   \frac{1}{2} \partial_t  \|\rho h \vt \|^2+  \rho_wC_{w1} \|\vt \|^2-(\div(\sigmat),\vt)=(\rhs,\vt),
\end{align}
{\GG with $\zeta_\text{min}:=\min_{(\xt,t) \in \Omega\times I} \zeta(\xt,t)$.}
Next we reformulate the stress tensor given in equation (\ref{stress})
\begin{align}
( -\div(\sigmat),\vt)=&( \sigmat, \nabla \vt)= \frac{1}{2}\big(\sigmat+\sigmat^T,\nabla\vt\big) 
    = \frac{1}{2}\big(\sigmat,\nabla\vt+\nabla\vt^T\big)=( \sigmat, \epsilont)\nonumber\\
=&\Big(2 \eta \epsilont + (\zeta-\eta) \tr(\epsilont) I-\frac{P}{2}I, \epsilont \Big)\nonumber\\
=&
\Big(\frac{1}{2} \zeta \epsilont +\frac{3}{4} \zeta \tr(\epsilont) I-\frac{P}{2}I, \epsilont \Big)\nonumber\\
=& \Big(\frac{1}{2} \zeta  \epsilont, \epsilont \Big)+\Big(\frac{3}{4} \zeta   \tr(\epsilont), \tr(\epsilont)\Big)\label{stress_est},
\end{align}
where we take into account that $\div(P)=0$ as $A$ and $h$ are constant.
We apply Korn's first inequality for homogeneous Dirichlet boundary values  \cite{Ciarlet2010} to obtain
\begin{align*}
 \|\epsilont \|^2 \geq c_k    \| \nabla \vt \|^2 ,
\end{align*}
where $c_k$ is the positive constant of Korn's inequality (\ref{Korn:first}). {This implies the following lower bound for (\ref{stress_est})}
{\small
\begin{align}\label{eq:eststress}
\|{\sqrt{2}^{-1}} \sqrt{\zeta} \epsilont \|^2+\frac{1}{2} \|\frac{\sqrt 3}{2} \sqrt{\zeta} \tr (\epsilont) \|^2 \geq  \|{\sqrt{2}^{-1}} \sqrt{\zeta}\epsilont \|^2 \geq  \frac{\zeta_\text{min}}{2}  \|\epsilont \|^2
\geq  c_k \frac{\zeta_\text{min}}{2} \|\nabla \vt \|^2.
\end{align}
}
{A combination of the $\vt$-dependent ocean drag (\ref{eq:oceandrag}) {and the estimate of the stress tensor (\ref{eq:eststress}) gives}
\begin{align}\label{eq:operator}
\Big(\rho_w C_{w1}\vt-\div(\sigmat),\vt\Big) \geq \rho_wC_{w1} \|\vt \|^2 +  c_k \frac{\zeta_\text{min}}{2} \|\nabla \vt \|^2.
\end{align}
}
{Finally we estimate the right hand side $\rhs$ of (\ref{eq:to_be_estimated}). With the inequalities of Cauchy-Schwarz and Young follows 
\begin{align}\label{eq:rhs}
( \rhs, \vt)  \leq  \frac{1}{2 \epsilon} \|\rhs\|^2+ \frac{\epsilon}{2}\| \vt \|^2 \leq {\G \frac{1}{2 \epsilon} \|\rhs\|^2
+c_p\frac{\epsilon}{2}  \| \nabla \vt \|^2},
\end{align}
where we have used in the last step Poincare's inequality with homogeneous Dirichlet boundaries \cite{evans10}. Here $c_p$ denotes the  {positive } constant of Poincare's inequality. {To move $\| \nabla \vt\|^2$ in (\ref{eq:rhs}) to the left hand side of equation (\ref{eq:operator}) we choose $\epsilon=c_k c_p^{-1} \frac{\zeta_\text{min}}{2}$}.
From the estimates (\ref{eq:operator}) and (\ref{eq:rhs}) follows for (\ref{eq:to_be_estimated})
\begin{align*}
 \partial_t  \|\rho h \vt \|^2+\frac{1}{2}\rho_wC_{w1} \|\vt \|^2 +  c_k \frac{\zeta_\text{min}}{4} \|\nabla \vt \|^2 \leq \frac{c_p}{c_k\zeta_\text{min}}  \|\rhs\|^2.
\end{align*}
We integrate over the time interval $I=[0,T]$ and the assertion follows
}
\begin{align*}
\|\rho h \vt(T) \|^2&- \|\rho h \vt(0) \|^2+ \int_{0}^T \{ \frac{1}{2}\rho_w C_{w1} \|\vt \|^2 +    c_k \frac{\zeta_\text{min}}{4} \|\nabla \vt \|^2 \} dt \\
& \leq \int_{0}^T \frac{c_p}{c_k\zeta_\text{min}}  \|\rhs\|^2 dt.
\end{align*}
\end{proof}
}

\begin{theorem}[$H^1$-estimate for EVP sea-ice momentum equation]\label{thoe2}
{\small
{Let a time interval [0,T], with $T>0$, be given. Let sea-ice concentration $A$ and mean sea-ice thickness $h$ be constant in time and space and{ $A$, $h>0$}.  Suppose
$\vt \in L^2([0,T],H^1_0(\Omega))$ 
is a solution of the sea-ice momentum equation (\ref{eq:mo}) with the EVP rheology given by (\ref{EVP}) and { $\tr(\sigmat)=0$}.}
Then $\vt$
satisfies the following $H^1$-estimate 
\begin{equation}
\begin{aligned}
 \rho_w C_{w1}  \|\vt(T) \|^2  &+c_k\frac{\zeta_\text{min}}{4} \|\nabla \vt(T)\|^2\\
&+\int_0^T\{\frac{1}{T_\text{EVP}} \rho_w C_{w1}\| \vt \|^2 +\frac{\zeta_\text{min}}{2 T_\text{EVP}}\|\nabla \vt \|^2  +\frac{T_\text{EVP}}{4} \rho_w C_{w1}\|\partial_t \vt \|^2
     \} dt\\
  & \leq    \frac{\zeta}{4}\|\epsilont(0)\|^2 + a \|\partial_t \vt(0)\|^2 +b \| \vt(0)\|^2
     + \int_0^T  c \|  \partial_t \rhs+\rhs \|^2 dt ,
 \end{aligned}
 \end{equation}
 with $a=\frac{\rho h T_\text{EVP}}{4}+ \frac{1}{2\gamma}\rho h$, $ b= \rho_w C_{w1}  +\frac{1}{T_\text{EVP}}\rho h + \frac{\gamma}{2}\rho h $, $ {c=  \frac{c_p \zeta_\text{min}}{2 c_k T_\text{EVP}} + \frac{T_\text{EVP}}{4 \gamma_2} }$. {\GG The minimal value of the viscosity is defined as $\zeta_\text{min}:=\min_{(\xt,t) \in \Omega\times I} \zeta(\xt,t)$} and }$\rhs= \rho_w C_{w1}\vt_w+ \rho_a C_a \|\vt_a\|\vt_a+\rho h f \te \times \vt_w$.
\end{theorem}
\begin{proof}
We reformulate the sea-ice momentum equation in the following form
\begin{align}\label{EVP:reduced}
\partial_t \rho h \vt +\rho_wC_{w1}\vt+ \rho h f_c \te \times \vt -\rhs&=\div(\sigmat).
\end{align}
We take the divergence of the elastic-viscous-plastic model approximation (\ref{EVP})  and apply that $tr(\sigmat)=0$. This gives
\begin{align}\label{EVP:div}
\partial_t \div(\sigmat)+ \frac{e^2}{2T_\text{EVP}}\div(\sigmat)=\frac{1}{T_\text{EVP}}\div(\zeta \epsilont).
\end{align}
Inserting equation (\ref{EVP:reduced}) into (\ref{EVP:div}) yields 
\begin{equation}
\begin{aligned}\label{EVP:divin}
  \rho h\partial_{t t}\vt+ \rho_w C_{w1} \partial_t\vt+\rho h f_c\te \times\partial_t \vt-\partial_t \rhs+&\\
  \frac{2}{T_\text{EVP}}\rho h \partial_t \vt +\frac{2}{T_\text{EVP}} \rho_w C_{w1}\vt+  \frac{2}{T_\text{EVP}} \rho h f_e \te \times \vt- \rhs&= \frac{1}{T_\text{EVP}}\div(\zeta \epsilont).
\end{aligned}
\end{equation}
Multiplying (\ref{EVP:divin}) with $\partial_t \vt$ and integrating  over the time interval $I=[0,T]$ results after division by $\frac{2}{T_\text{EVP}}$ in
\begin{equation}\label{EVP:vt}
 \begin{aligned}
\int_I \Big\{ \, \medskip  \frac{\rho h T_\text{EVP}}{4} \partial_t\|\partial_t \vt\|^2+\frac{T_\text{EVP}}{2}\rho_w C_{w1}\|\partial_t \vt \|^2 -&\\
\frac{T_\text{EVP}}{2}\ \partial_t \rhs\partial_t \vt
    +\rho h \|\partial_t \vt\|^2+ \frac{\rho_w}{2} C_{w1} \partial_t \|\vt\|^2+&\\
    \rho h f_c \te \times \vt \partial_t \vt   - \frac{T_\text{EVP}}{2}\rhs \partial_t \vt  &\Big\} dt  =  \int_I \{-\frac{\zeta}{4}\partial_t \|\epsilont\|^2 \}dt.
 \end{aligned}
 \end{equation}
Multiplication of (\ref{EVP:divin}) with  $\vt$ and using that the multiplication with the Coriolis term is anticommutative yields
 \begin{equation}\label{EVP:v}
 \begin{aligned}
 \rho h \partial_t \vt(T)\vt(T) -\rho h \partial_t \vt(0)\vt(0)+&\\
\int_I \Big\{ -\rho h \|\partial_{t}\vt \|^2 +\frac{\rho_w}{2} C_{w1} \partial_t \| \vt \|^2  - \rho h f_c \te \times \vt \partial_t \vt-&\\  \partial_t \rhs \vt+
  \frac{1}{T_\text{EVP}}\rho h \partial_t\| \vt\|^2 +\frac{2}{T_\text{EVP}}\rho_wC_{w1} \| \vt \|^2  -\rhs \vt &\Big\} dt =\int_I \{-\frac{\zeta}{T_\text{EVP}} \|\epsilont \|^2\} dt.
 \end{aligned}
 \end{equation}
 Adding (\ref{EVP:vt}) and (\ref{EVP:v}) results in
\begin{equation}
 \begin{aligned}
 \rho h \partial_t \vt(T)\vt(T) -\rho h \partial_t \vt(0)\vt(0)+   
 \int_I \Big \{ 
 \frac{T_\text{EVP}}{2}\rho_w C_{w1}\|\partial_t \vt \|^2 +&\\
 \frac{\rho h T_\text{EVP}}{4} \partial_t\|\partial_t \vt\|^2
 -
 \frac{T_\text{EVP}}{2} (\partial_t \rhs +\rhs)\partial_t \vt+
 \rho_w C_{w1} \partial_t \|\vt\|^2 +&\\
  \frac{1}{T_\text{EVP}}\rho h \partial_t\| \vt\|^2 +
  \frac{2}{T_\text{EVP}}\rho_wC_{w1} \| \vt \|^2  -( \partial_t \rhs+\rhs) \vt 
    \Big \} dt=& \int_I \Big \{-\frac{\zeta}{4}\partial_t \|\epsilont\|^2\\ &-\frac{\zeta}{T_\text{EVP}}\|\epsilont \|^2 \Big\} dt  
 \end{aligned}
 \end{equation}
 We integrate over time and apply Young's inequality to get
 \begin{equation}
 \begin{aligned}
 -&\frac{\gamma}{2}\rho h\| \partial_t \vt(T)\|^2- \frac{1}{2\gamma}\rho h \|\vt(T)\|^2  -\frac{1}{2\gamma}\rho h \|\partial_t \vt(0)\|^2- \frac{\gamma}{2}\rho h \|\vt(0)\|^2\\
 &+ \frac{\rho h T_\text{EVP}}{4} \|\partial_t \vt(T)\|^2+ \rho_w C_{w1}  \|\vt\|^2 (T)  +\frac{1}{T_\text{EVP}}\rho h \| \vt(T)\|^2 +\frac{\zeta}{4} \|\epsilont(T)\|^2\\    
&+\int_I\Big\{\frac{2}{T_\text{EVP}}\rho_wC_{w1} \| \vt \|^2  - { (\partial_t \rhs+\rhs)\vt} +\frac{\zeta}{T_\text{EVP}}\|\epsilont \|^2\\ &+\frac{T_\text{EVP}}{2}\rho_w C_{w1}\|\partial_t \vt \|^2 -\frac{T_\text{EVP}}{4 \gamma_2} \|(\partial_t \rhs +\rhs)\|^2 -\frac{\gamma_2 T_\text{EVP}}{4}\|\partial_t \vt \|^2
     \Big \}dt\\
     \leq
     &\frac{\zeta}{4}\|\epsilont(0)\|^2 + \frac{\rho h T_\text{EVP}}{4} \|\partial_t \vt(0)\|^2 + \rho_w C_{w1}  \|\vt(0)\|^2 +\frac{1}{T_\text{EVP}}\rho h \| \vt(0)\|^2,
 \end{aligned}
 \end{equation}
with $\gamma=\frac{T_\text{EVP}}{2}$ and $\gamma_2=\rho_wC_{w1}$. 
{Applying Korn's inequality (\ref{Korn:first}) to the strain rate tensor and using Poincare's inequality for $(\partial_t \rhs+ \rhs)\vt$ gives final estimate}
 \begin{equation}
 \begin{aligned}
&\rho_w C_{w1}  \|\vt(T) \|^2  +c_k\frac{\zeta_\text{min}}{4} \|\nabla \vt(T)\|^2\\
&+\int_I\Big \{\frac{2}{T_\text{EVP}}\rho_wC_{w1} \| \vt \|^2 { +c_k \frac{\zeta_\text{min}}{2 T_\text{EVP}}}\|\nabla \vt \|^2  +\frac{T_\text{EVP}}{4}\rho_w C_{w1}\|\partial_t \vt \|^2
     \Big \} dt\\
     \leq& 
     \frac{\zeta}{4}\|\epsilont(0)\|^2 + (\frac{\rho h T_\text{EVP}}{4} +\frac{1}{2\gamma}\rho h ) \|\partial_t \vt(0)\|^2 \\
     &+ (\rho_w C_{w1}  +\frac{1}{T_\text{EVP}}\rho h + \frac{\gamma}{2}\rho h ) \| \vt(0)\|^2
     + \int_I  (\frac{1}{2\gamma_0} + \frac{T_\text{EVP}}{4 \gamma_2} )\|\partial_t \rhs+\rhs\|^2\, dt,
 \end{aligned}
 \end{equation}
 {with $ \gamma_0=\frac{c_k\zeta_{min}}{c_p T_\text{EVP}}$ and {\GG  $\zeta_\text{min}:=\min_{(\xt,t) \in \Omega\times I} \zeta(\xt,t)$}.}
\end{proof}

\end{appendix}


\end{document}